\pdfoutput=1
\documentclass[a4paper,11pt]{amsart}

\usepackage{a4wide}
\usepackage{amsmath, amsfonts, amscd, amssymb, amsthm}
\usepackage{tikz}
\usepackage[UKenglish]{babel}
\usetikzlibrary{matrix,arrows,calc}
\usepackage{xspace} 
\usepackage[bottom]{footmisc}

\usepackage{aliascnt}

\theoremstyle{plain}
\newtheorem{thm}{Theorem}[section]

\newaliascnt{lem}{thm}
\newtheorem{lem}[lem]{Lemma}
\aliascntresetthe{lem}

\newaliascnt{prop}{thm}
\newtheorem{prop}[prop]{Proposition}
\aliascntresetthe{prop}               

\newaliascnt{cor}{thm}
\newtheorem{cor}[cor]{Corollary}
\aliascntresetthe{cor}

\theoremstyle{definition}
\newaliascnt{defn}{thm}
\newtheorem{defn}[defn]{Definition}
\aliascntresetthe{defn}

\newaliascnt{rem}{thm}
\newtheorem{rem}[rem]{Remark}
\aliascntresetthe{rem}

\newaliascnt{exa}{thm}
\newtheorem{exa}[exa]{Example}
\aliascntresetthe{exa}

\newcommand{\N}{\mathbb{N}}
\newcommand{\Z}{\mathbb{Z}}
\newcommand{\R}{\mathbb{R}}
\newcommand{\C}{\mathbb{C}}

\newcommand{\bi}{\begin{itemize}}
\newcommand{\ei}{\end{itemize}}
\newcommand{\be}{\begin{equation*}}
\newcommand{\ee}{\end{equation*}}
\newcommand{\bp}[1][]{\begin{proof}[Proof#1.]}				
\newcommand{\ep}{\end{proof}}
\newcommand{\wt}{\widetilde}

\newcommand{\mf}{\mathfrak}

\newcommand{\mc}{\mathcal}
\newcommand{\symp}{symplectic structure\xspace}

\newcommand{\blog}{log-\symp}

\newcommand{\nsymp}{near-\symp}

\newcommand{\acs}{almost-complex structure\xspace}

\newcommand{\gcs}{generalized complex structure\xspace}
\newcommand{\Gcs}{Generalized complex structure\xspace}
\newcommand{\sgcs}{stable \gcs}
\newcommand{\Sgcs}{Stable \gcs}

\newcommand{\alf}{achiral Lefschetz fibration\xspace}

\newcommand{\lf}{Lefschetz fibration\xspace}

\newcommand{\btwo}[1][]{{^b}{#1}} 
\renewcommand{\btwo}[1][]{#1} 

\newcommand{\etwo}[1][]{{^e}{#1}} 
\renewcommand{\etwo}[1][]{#1} 

\newcommand{\placeholder}{boundary \lf}

\numberwithin{equation}{section}							
\begin{document}

\author{Gil R.\ Cavalcanti}
\address{Department of Mathematics, Utrecht University, 3508 TA Utrecht, The Netherlands}
\email{g.r.cavalcanti@uu.nl}

\author{Ralph L.\ Klaasse}
\address{Department of Mathematics, Utrecht University, 3508 TA Utrecht, The Netherlands}
\email{r.l.klaasse@uu.nl}

\date{\today}
\thanks{This project was supported by a VIDI grant from NWO, the Dutch science foundation.}

\begin{abstract} A \gcs is called stable if its defining anticanonical section vanishes transversally, on a codimension-two submanifold. Alternatively, it is a zero elliptic residue \symp in the elliptic tangent bundle associated to this submanifold. We develop Gompf--Thurston symplectic techniques adapted to Lie algebroids, and use these to construct \sgcs{}s out of log-symplectic structures. In particular we introduce the notion of a \placeholder for this purpose and describe how they can be obtained from genus one Lefschetz fibrations over the disk.
\end{abstract}
\title{Fibrations and stable generalized complex structures}
\maketitle

\vspace{-2em}
\tableofcontents
\vspace{-3.5em}
\section{Introduction}
\label{sec:introduction}
Generalized complex structures \cite{Gualtieri04, Gualtieri11, Hitchin03} capture aspects of both symplectic and complex geometry. While slightly misleading, a generalized complex structure can be seen as a Poisson structure together with a suitably compatible complex structure normal to the (singular) symplectic leaves. In general the symplectic leaves may have varying dimension, leading to the notion of type change, where complex and symplectic behavior is mixed.

In this paper we are interested in a class of generalized complex manifolds for which the type need not be constantly equal to zero, but differs from zero only in the mildest way possible. These are called \sgcs{}s \cite{CavalcantiGualtieri15,Goto15}, and are those \gcs{}s whose defining anticanonical section vanishes transversally. Whenever a \gcs is of type zero, it is isomorphic after a $B$-field transformation to a symplectic structure. In this sense a \sgcs is a symplectic-like structure that fails to be symplectic on at most a codimension-two submanifold.

In \cite{CavalcantiGualtieri09, GotoHayano16, Torres12, TorresYazinski14}, many examples of \sgcs{}s on four-manifolds were constructed, in particular on manifolds without a symplectic nor complex structure. A thorough study of stable generalized complex manifolds was initiated in \cite{CavalcantiGualtieri15}. \Sgcs{}s can be seen as the generalized geometric analogue of a \blog, which is a type of mildly degenerate Poisson structure that has recently received a lot of attention \cite{Cavalcanti13, CavalcantiKlaasse16, GualtieriLi14, GuilleminMirandaPires14, MarcutOsornoTorres14, MarcutOsornoTorres14two}.

In \cite{CavalcantiGualtieri15} it was noted that \sgcs{}s can alternatively be viewed as symplectic forms on a Lie algebroid, the elliptic tangent bundle, constructed out of the anticanonical section and its zero set (they are further of zero elliptic residue: see Section \ref{sec:gcsstablegcs} for details.). This symplectic viewpoint allows one to apply symplectic techniques to the study of \sgcs{}s. In this paper we focus on constructing \sgcs{}s on total spaces of fibration-like structures.

More precisely, in this paper we extend Gompf--Thurston symplectic techniques to arbitrary Lie algebroids. We introduce Lie algebroid fibrations and Lie algebroid Lefschetz fibrations, and give criteria when these can be equipped with compatible Lie algebroid symplectic structures. We discuss morphisms between elliptic and log-tangent bundles, which are Lie algebroids whose symplectic structures describe \sgcs{}s and \blog{}s respectively. We then introduce a class of maps called boundary maps. These induce morphisms between the elliptic and log tangent bundles, and we use them to construct \sgcs{}s out of \blog{}s.

The following is the main result about Lie algebroid Lefschetz fibrations and Lie algebroid symplectic structures. A more precise version exists (\autoref{thm:thurstonlfibration}), as well as an analogous result for Lie algebroid fibrations (\autoref{thm:thurstonfibration}).
\begin{thm} Let $(\varphi,f)\colon \mc{A}^4_{X} \to \mc{A}^2_{\Sigma}$ be a Lie algebroid \lf with connected fibers. Assume that $\mc{A}_\Sigma$ admits a symplectic structure and there exists a closed $\mc{A}_X$-two-form $\eta$ such that $\eta|_{\ker \varphi}$ is nondegenerate. Then $X$ admits an $\mc{A}_X$-\symp.
\end{thm}
A specific type of boundary map we call a \placeholder{} is shown to induce a Lie algebroid Lefschetz fibration between the elliptic and log-tangent bundle. This concept extends and formalizes the generalized Lefschetz fibrations of \cite{CavalcantiGualtieri09}. Using this correspondence we are able to prove the following result (\autoref{thm:lfibrationsgcs}) which motivated this work.
\begin{thm} Let $f\colon (X^4,D) \to (\Sigma,Z)$ be a boundary Lefschetz fibration such that $(\Sigma,Z)$ carries a \blog. Assume that ${[F] \neq 0 \in H_2(X \setminus D;\R)}$ for the generic fiber $F$. Then $(X,D)$ admits a \sgcs.
\end{thm}
A similar result is true in arbitrary dimension, using boundary fibrations instead (\autoref{thm:fibrationsgcs}). This result validates the definition of a boundary Lefschetz fibration as being the type of fibration-like map which should be linked with \sgcs{}s. This extends the original relationship between Lefschetz fibrations and symplectic structures due to Gompf \cite{GompfStipsicz99}, that between broken \lf{}s and \nsymp{}s \cite{AurouxDonaldsonKatzarkov05, Baykur09}, and between \alf{}s and \blog{}s \cite{Cavalcanti13,CavalcantiKlaasse16}.

To produce concrete examples, one has to construct \placeholder{}s on explicit four-manifolds. This can be done using genus one Lefschetz fibrations over punctured surfaces with boundary monodromy given by powers of Dehn twists using a process we call completion (see \autoref{cor:lfcompletion}). In particular, we have the following corollary (\autoref{thm:g1lfsgcs}).
\begin{thm} Let $f\colon X \to D^2$ be a genus one Lefschetz fibration over the disk whose monodromy around the boundary is a power of a Dehn twist. Then all possible completions $\wt{f}\colon (\wt{X}, D) \to (D^2, \partial D^2)$ admit a \sgcs.
\end{thm}
Together with Stefan Behrens we are working on classifying boundary Lefschetz fibrations over the disk.

All manifolds in this paper will be compact and orientable, unless specifically stated otherwise. Note however that they are not necessarily closed. Throughout we will identify de Rham cohomology with singular cohomology with $\R$-coefficients.
\subsection*{Organization of the paper}
In Section \ref{sec:gcsstablegcs} we recall the definition of a \sgcs and introduce the language of divisors. We further define Poisson structures of divisor-type and recall the definition of a \blog.

In Section \ref{sec:lalgebroids} we recall the definition of a Lie algebroid, and introduce the two relevant Lie algebroids (the elliptic and log tangent bundle) constructed out of divisors. We further define Lie algebroid fibrations and Lie algebroid Lefschetz fibrations, and discuss Lie algebroid symplectic structures and their Poisson counterparts.

In Section \ref{sec:gompfthurston} we discuss the Lie algebroid version of the Gompf--Thurston argument to construct Lie algebroid symplectic structures (\autoref{thm:lathurstontrick}). Moreover, we prove results for Lie algebroid fibrations (\autoref{thm:thurstonfibration}) and Lie algebroid Lefschetz fibrations (\autoref{thm:thurstonlfibration}).

In Section \ref{sec:boundarymaps} we define the normal Hessian of a map and introduce boundary maps. We define boundary fibrations and boundary Lefschetz fibrations and prove several normal form results for boundary maps (\autoref{prop:bdrymaplocalcoord}, \autoref{prop:fibratingpointwise}, and \autoref{prop:fibratingsemiglobal}).

In Section \ref{sec:constrbdrlfs} we introduce a standard boundary fibration (\autoref{prop:standardbdryfibration}) and use this to obtain boundary Lefschetz fibrations out of genus one Lefschetz fibrations over a punctured surface using monodromy data (\autoref{cor:lfcompletion}) via a completion process.

In Section \ref{sec:constrsgcs} we then prove our main results: \autoref{thm:lfibrationsgcs}, that four-dimensional boundary Lefschetz fibrations give rise to \sgcs{}s; and \autoref{thm:fibrationsgcs}, that the same holds for boundary fibrations. Moreover, we show compatibility with existing fibrations over $T^2$ and $S^1$ can be achieved in the compact case (\autoref{cor:compatiblefibrations}).

In Section \ref{sec:applications} we give examples of stable generalized complex manifolds constructed using our methods. In particular, we recover the examples $m \C P^2 \# n \overline{\C P}^2$ for $m$ odd from \cite{CavalcantiGualtieri09}.
\subsection*{Acknowledgements}
Apart from financial support mentioned on the title page, the authors would like to thank Selman Akbulut, Stefan Behrens, Robert Gompf and Andr\'as Stipsicz for useful discussions.
\section{Stable generalized complex and log-symplectic structures}
\label{sec:gcsstablegcs}
In this section we recall the notion of a \sgcs as defined in \cite{CavalcantiGualtieri15}, and that of \blog{}s \cite{GuilleminMirandaPires14}. Moreover, we introduce the language of divisors to study them. For a more comprehensive account, see also \cite{Klaasse17}. We start by recalling the definition of a \gcs, for which a general reference is \cite{Gualtieri11}.
\subsection{\Gcs{}s}
Let $X$ be a $2n$-dimensional manifold equipped with a closed three-form $H \in \Omega^3_{\rm cl}(X)$. Recall that the \emph{double tangent bundle} $\mathbb{T}X := TX \oplus T^*X$ is a Courant algebroid whose anchor is the projection $p\colon \mathbb{T}X \to TX$.  It carries a natural pairing $\langle V + \xi, W + \eta\rangle = \frac12(\eta(V) + \xi(W))$ of split signature and an $H$-twisted Courant bracket $[[V + \xi, W + \eta]]_H = [V, W] + \mc{L}_V \eta - \iota_W d\xi + \iota_V \iota_W H$ for $V, W \in \Gamma(TX)$ and $\xi, \eta \in \Gamma(T^*X)$.
\begin{defn} A \emph{\gcs} on $(X,H)$ is a complex structure $\mc{J}$ on $\mathbb{T}X$ that is orthogonal with respect to $\langle\cdot,\cdot\rangle$, and whose $+i$-eigenbundle is involutive under $[[\cdot,\cdot]]_H$.
\end{defn}
There is an alternative definition of a \gcs using spinors. To state it, recall that sections $v = V + \xi \in \Gamma(\mathbb{T}X)$ of the double tangent bundle act on differential forms via Clifford multiplication, given by $v \cdot \rho = \iota_V \rho + \xi \wedge \rho$ for $\rho \in \Omega^\bullet(X)$.
\begin{defn} A \emph{\gcs} on $(X,H)$ is given by a complex line bundle $K_{\mc{J}} \subset \wedge^\bullet T^*_\C X$ pointwise generated by a differential form $\rho = e^{B + i \omega} \wedge \Omega$ with $\Omega$ a decomposable complex $k$-form, satisfying $\Omega \wedge \overline{\Omega} \wedge \omega^{n-k} \neq 0$, and such that $d \rho + H \wedge \rho = v \cdot \rho$ for any local section $\rho \in \Gamma(K_{\mc{J}})$ and some $v \in \Gamma(\mathbb{T}X)$.
\end{defn}
Both definitions are related using that $K_\mc{J} = {\rm Ann}(E_\mc{J})$ is the annihilator under the Clifford action of $E_\mc{J}$, the $+i$-eigenbundle of $\mc{J}$. The bundle $K_{\mc{J}}$ is called the \emph{canonical bundle} of $\mc{J}$. For later use, we further introduce the analogue of a Calabi--Yau manifold in generalized geometry. Denote by $d^H = d + H \wedge$ the $H$-twisted de Rham differential.
\begin{defn} A generalized complex structure $\mc{J}$ on $(X,H)$ is \emph{generalized Calabi--Yau} if its canonical bundle $K_\mc{J}$ is determined by a global nowhere vanishing $d^H$-closed form.
\end{defn}
\begin{exa}\label{exa:gcss} The following provide examples of \gcs{}s on $(X,0)$.
\bi
	\item Let $\omega$ be a symplectic structure on $X$. Then $K_{\mc{J}_\omega} := \langle e^{i \omega} \rangle$ defines a \gcs $\mc{J}_\omega$.
	\item Let $J$ be a complex structure on $X$ with canonical bundle $K_J = \wedge^{n,0} T^* X$. Then $K_{\mc{J}_J} := K_J$ defines a \gcs $\mc{J}_J$.
	\item Let $P \in \Gamma(\wedge^{2,0} TX)$ a holomorphic Poisson structure with respect to a complex structure $J$. Then $K_{\mc{J}_{P,J}} := e^{P} K_J$ defines a \gcs $\mc{J}_{P,J}$.
\ei
The automorphisms $\mc{J}_\omega$, $\mc{J}_J$ and $\mc{J}_{P,J}$ are given by, with $\pi = {\rm Im}(P)$:
\begin{align*}
	\mc{J}_\omega = \begin{pmatrix} 0 & -\omega^{-1} \\ \omega & 0\end{pmatrix},  \qquad\qquad &	\mc{J}_J = \begin{pmatrix} -J & 0 \\ 0 & J^* \end{pmatrix}, \qquad \qquad \mc{J}_{P,J} = \begin{pmatrix} -J & \pi \\ 0 & J^* \end{pmatrix}.
\end{align*}	
\end{exa}
We next introduce the type of a \gcs$\mc{J}$, which colloquially provides a measure for how many complex directions there are. The type is an integer-valued upper semicontinuous function on $X$ whose parity is locally constant.
\begin{defn} Let $\mc{J}$ be a \gcs on $X$. The \emph{type} of $\mc{J}$ is a map ${\rm type}(\mc{J})\colon X \to \Z$ whose value at a point $x \in X$ is the integer $k$ above, the degree of $\Omega$. The \emph{type change locus} $D_{\mc{J}}$ of $\mc{J}$ is the subset of $X$ where ${\rm type}(\mc{J})$ is not locally constant.
\end{defn}
The complement of $D_{\mc{J}}$ is an open dense set where the type is minimal. Using the type, \gcs{}s are seen to interpolate between symplectic and complex structures. At points where ${\rm type}(\mc{J}) = 0$, the \gcs is equivalent to a \symp, in that it is equivalent to the \gcs $\mc{J}_\omega$ of a \symp under a \emph{$B$-field transformation}. These are automorphisms $e^B\colon \mathbb{T}X \to \mathbb{T}X$ given by $e^B\colon V + \xi \mapsto V + \xi + \iota_V B$ for $B \in \Omega^2_{\rm cl}(X)$.

Any \gcs $\mc{J}$ determines a Poisson structure $\pi_\mc{J}$ as the composition $\pi_{\mc{J}}^\sharp := p \circ \mc{J}|_{T^*X}$. The type of $\mc{J}$ is related to the rank of $\pi_{\mc{J}}$ through the formula ${\rm rank}(\pi_{\mc{J}}) = 2n - 2 \, {\rm type}(\mc{J})$. Using $\pi_{\mc{J}}$ one can view a \gcs $\mc{J}$ as a foliation on $X$ with symplectic leaves, and a suitably compatible complex structure transverse to the leaves.
\subsection{\Sgcs{}s}
\Gcs{}s which are stable were introduced in \cite{CavalcantiGualtieri15,Goto15}. Their defining property is a natural condition and since \sgcs{}s are not far from being symplectic, one can use symplectic techniques to study them.

Let $\mc{J}$ be a \gcs on $(X,H)$. The anticanonical bundle $K_{\mc{J}}^*$ has a section $s \in \Gamma(K_{\mc{J}}^*)$, given by $s(\rho) := \rho_0$ for $\rho \in \Gamma(K_{\mc{J}})$, with $\rho_0$ the degree-zero part of $\rho$.
\begin{defn} A \gcs $\mc{J}$ on $(X,H)$ is \emph{stable} if $s$ is transverse to the zero section in $K_{\mc{J}}^*$. The set $D_{\mc{J}} := s^{-1}(0)$ is a codimension-two smooth submanifold of $X$ called the \emph{anticanonical divisor} of $\mc{J}$.
\end{defn}
Outside of $D_{\mc{J}}$, the section $s$ is nonvanishing hence the type of $\mc{J}$ is equal to zero, while over $D_{\mc{J}}$ it is equal to two. Consequently, \sgcs{} can be seen as \gcs{}s which are close to being symplectic.

\begin{exa} Consider $(\C^2,0)$ with holomorphic Poisson structure $\pi = z \partial_z \wedge \partial_w$. This gives a \sgcs with $K_{\mc{J}} = \langle z + dz \wedge dw \rangle$ and $D_\mc{J} = \{z = 0\}$.
\end{exa}
By \autoref{exa:gcss}, any holomorphic Poisson structure $(J,P)$ defines a \gcs $\mc{J}_{P,J}$ on $(X,0)$ with $K_{\mc{J}_{P,J}} = e^P K_J$. Thus $K_{\mc{J}_{P,J}}$ is locally generated by $e^P \Omega$, with $\Omega$ a local trivialization of $K_J$. This immediate proves the following.
\begin{prop}[{\cite[Example 2.17]{CavalcantiGualtieri15}}] Let $(X,0,J,P)$ be a holomorphic Poisson structure on a complex $2m$-dimensional manifold. Then $\mc{J}_{P,J}$ is a \sgcs if and only if the Pfaffian $\wedge^m P$ is transverse to the zero section in $\Gamma(\wedge^{2m,0} TX)$.
\end{prop}
Any \sgcs $\mc{J}$ is locally equivalent around points in $D_\mc{J}$ to $\langle e^{i \omega_0}(z + dz \wedge dw) \rangle$ on $\C^2 \times \R^{2n-4}$, with $\omega_0$ the standard symplectic form on $\R^{2n-4}$, and $D_{\mc{J}} = \{z = 0\}$. On $D_\mc{J}$ there is an induced type-$1$ generalized Calabi--Yau structure.

Any compact type-1 generalized Calabi--Yau manifold, such as $D_\mc{J}$, fibers over the torus $T^2$ \cite{BaileyCavalcantiGualtieri16}. Moreover, the semilocal form of a \sgcs around its type change locus is given by its \emph{linearization} along $D_{\mc{J}}$, which is the \sgcs naturally present on the normal bundle to this type-1 generalized Calabi--Yau manifold. We will not elaborate on this further and instead refer to \cite{BaileyCavalcantiGualtieri16, CavalcantiGualtieri15}.

The anticanonical bundle $K^*_{\mc{J}}$ of a stable generalized complex manifold together with its natural section $s$ are a particular example of a \emph{divisor}, which we will introduce shortly. As the theory of divisors permeates much of this work, we now turn to developing this concept.

\subsection{Real and complex divisors}
\label{sec:divisors}
In this section we introduce real and complex divisors, which are extensions to the smooth setting of the notion of a divisor used in complex geometry. We further discuss the relation between the ideals they define, and define morphisms between divisors. See also \cite{CavalcantiGualtieri15} and \cite{VanderLeerDuran16}.
\begin{defn} Let $X$ be a manifold. A \emph{divisor} on $X$ is a pair $(U,\sigma)$ where $U \to X$ is a real/complex line bundle and $\sigma \in \Gamma(U)$ is a section whose zero set is nowhere dense.
\end{defn}
In this paper we will mostly focus on real divisors. As such, we will often drop the prefix `real', while explicitly stating when divisors are instead complex. Examples of divisors will be discussed in upcoming subsections.

Given a divisor $(U,\sigma)$, viewing $\sigma\colon \Gamma(U^*) \to C^\infty(X)$ specifies a locally principal ideal $I_\sigma := \sigma(\Gamma(U^*))$ which is locally generated by a function with nowhere dense zero set. Letting $\alpha$ be a local trivialization of $U^*$, we have $\alpha(\sigma) = g$ for some local function $g$. Then locally $I_\sigma = \langle \alpha(\sigma) \rangle = \langle g \rangle$.
Conversely, out of any such ideal we can construct a divisor, which recovers the ideal via this evaluation process. This extends the correspondence between divisors and holomorphic line bundles in complex geometry.
\begin{prop}[\cite{VanderLeerDuran16}]\label{prop:locprincideal} Let $I$ be a locally principal ideal on $X$ generated by functions with nowhere dense zero set. Then there exists a divisor $(U_I,\sigma)$ on $X$ such that $I_\sigma = I$.
\end{prop}
\bp Let $\{U_\alpha\}$ be an open cover of $X$ and $f_\alpha \in I(U_\alpha)$ be generators. Then on $U_\alpha \cap U_\beta$ we have $f_\alpha = g_{\alpha \beta} f_\beta$ with $g_{\alpha \beta} \in C^\infty(U_\alpha \cap U_\beta)$. Since $f_\alpha = g_{\alpha\beta} g_{\beta \alpha} f_\alpha$, and $f_\alpha$ has nowhere dense support, we see that each $g_{\alpha\beta}$ is a nonvanishing function and $g_{\beta\alpha} = g_{\alpha\beta}^{-1}$. Similarly the identity $f_\alpha = g_{\alpha \beta} g_{\beta \gamma} g_{\gamma \alpha} f_\alpha$ on $U_\alpha \cap U_\beta \cap U_\gamma$ implies that the functions $g_{\alpha \beta}$ satisfy the cocycle condition. We conclude that $\{(U_{\alpha \beta}, g_{\alpha \beta})\}$ defines a line bundle $U_I$ on $X$. Further, setting $\sigma|_{U_\alpha} = f_\alpha$ on $U_\alpha$ specifies a section $\sigma$ of this bundle with the desired properties.
\ep
The section $\sigma$ constructed in the proof of the above proposition is unique up to multiplication by a smooth nonvanishing function. We next define morphisms between divisors in terms of the ideals they give rise to. Denote by $f^*I \subseteq C^\infty(X)$ the ideal generated by the pullback of an ideal $I \subseteq C^\infty(Y)$ along a map $f\colon X \to Y$.
\begin{defn} Let $(U_X,\sigma_X)$ be a divisor on $X$ and $(U_Y,\sigma_Y)$ a divisor on $Y$. A map $f\colon X \to Y$ is a \emph{morphism of divisors} if $f^* I_{\sigma_Y} = I_{\sigma_X}$.
\end{defn}
Equivalently, one can require that $(U_X, \sigma_X) = (f^* U_Y, g f^* \sigma_Y)$ for some $g\in C^{\infty}(X;\R^*)$. Two divisors are \emph{diffeomorphic} (denoted using $\cong$) if there exists a morphism of divisors between them which is in addition a diffeomorphism. Two divisors on a fixed manifold $X$ are \emph{isomorphic} (denoted using $=$) if the identity map on $X$ is a morphism of divisors.
\subsubsection{Log divisors}
The simplest type of divisor is that of a log divisor. These are divisors whose zero set will be a smooth codimension-one submanifold, or \emph{hypersurface}.
\begin{defn}\label{defn:logdivisor} A \emph{log divisor} is a divisor $(L,s)$ whose zero set $Z_s = s^{-1}(0)$ is a smooth hypersurface along which $s$ is transverse to the zero section.
\end{defn}
It follows from the definition that the first Stiefel--Whitney class of $L$ is the $\Z_2$-Poincar\'e dual of $Z$. Further, the intrinsic derivative $ds|_Z\colon NZ \to L|_Z$ is an isomorphism. 

There is also a complex analogue of a (real) log divisor.
\begin{defn} A \emph{complex log divisor} is a complex divisor $(U,\sigma)$ whose zero set ${D_\sigma = \sigma^{-1}(0)}$ is a smooth codimension-two submanifold with $\sigma$ transverse to the zero section.
\end{defn}
Any complex log divisor $(U,\sigma)$ has a complex conjugate $(\overline{U}, \overline{\sigma})$ with the same zero locus.
\begin{exa} By definition the anticanonical bundle $K_\mc{J}^*$ of a \sgcs together with its natural section $s$ is an example of a complex log divisor.
\end{exa}
Let $(L,s)$ be a log divisor. The associated ideal $I_Z := I_s$ is exactly the vanishing ideal of the hypersurface $Z$. Any hypersurface naturally gives rise to a unique log divisor.
\begin{prop}\label{prop:hypersurfacedivisor} Let $X$ be a manifold and $Z \subset X$ a hypersurface. Then $Z$ carries a unique log divisor structure, i.e.\ there exists a unique log divisor $(L,s)$ with $s^{-1}(0) = Z$.
\end{prop}
\bp Apply \autoref{prop:locprincideal} to the vanishing ideal $I_Z$, giving a divisor $(U_{I_Z}, \sigma) =: (L,s)$. The section $s$ vanishes transversely along $Z$ as it is equal to a local defining function for $Z$ in any trivializing open $U_\alpha$ of $L$ containing $Z$. We conclude that $(L,s)$ is a log divisor.
\ep
Because of this result we often identify a log divisor with the associated submanifold $Z$.

\subsubsection{Elliptic divisors}
More directly relevant to our study of \sgcs{}s is the notion of an elliptic divisor.
\begin{defn} An \emph{elliptic divisor} is a divisor $(R,q)$ whose zero set $D_q = q^{-1}(0)$ is a smooth codimension-two critical submanifold of $q$ along which its normal Hessian is definite.
\end{defn}
We denote an elliptic divisor by $|D|$. The normal Hessian ${\rm Hess}(q) \in \Gamma(D; {\rm Sym}^2 N^*D \otimes R)$ of $q$ is the leading term of its Taylor expansion. There is more information available in the elliptic divisor than just the zero set of $q$. The ideal $I_D := I_q$ is called an \emph{elliptic ideal}, and both $R$ and $q$ (up to a nonzero smooth function) can be recovered from $I_D$ by \autoref{prop:locprincideal}.
\begin{exa} Let $D = (U,\sigma)$ be a complex log divisor. Then $((U\otimes \overline{U})_\R, \sigma \otimes \overline{\sigma})$ is an elliptic divisor $|D|$. Using \autoref{prop:ellmorsebott} below and the factorization $x^2 + y^2 = (x + iy)(x - iy) = w \overline{w}$, any elliptic divisor with coorientable zero set arises from a complex log divisor in this way, with $(U,\sigma)$ being determined upto diffeomorphism by the choice of coorientation.
\end{exa}
Note that $I_D$ is not the vanishing ideal of $D$, but instead is locally generated by an even index Morse--Bott function in coordinates normal to $D$, as we now explain.
\begin{defn} Let $g \in C^\infty(X)$ be given. A compact connected submanifold $S \subset X$ is a \emph{nondegenerate critical submanifold} of $g$ if $S \subset {\rm Crit}(g)$ and $\ker {\rm Hess}(g) = T_p S$ for all $p \in S$. If ${\rm Crit}(g)$ consists of nondegenerate critical submanifolds, then $g$ is a \emph{Morse--Bott function}.
\end{defn}
Let $g \in C^\infty(X)$ and take $S$ a nondegenerate critical submanifold of $g$. Consider the exact sequence $0 \to TS \to TX|_S \to NS \to 0$. For $p \in S$ we have ${\rm Hess}(g)(p) \in {\rm Sym}^2 T_p^* X$, and this vanishes when restricted to $T_p S$. But then ${\rm Hess}(g)(p) \in {\rm Sym}^2 N_p^* S$, giving a nondegenerate bilinear form $Q_g \in \Gamma(S; {\rm Sym}^2 N^* S)$. The semi-global version of the Morse--Bott lemma says that $g$ is diffeomorphic to this quadratic approximation in a tubular neighbourhood of $S$.
\begin{lem}[{\cite[Proposition 2.6.2]{Nicolaescu11}}]\label{lem:morsebott} Let $g \in C^\infty(X)$ be a Morse--Bott function and $S$ a nondegenerate critical submanifold of $g$. Then there exists a neighbourhood $U$ of the zero section $S \subset NS$ and an open embedding $\Phi\colon U \to X$ such that $\Phi|_S = {\rm id}_S$ and $\Phi^* g = Q_g$.
\end{lem}
Let $|D| = (R,q)$ be an elliptic divisor. Then $R$ is orientable by $q$ as it is a trivialization away from a codimension-two submanifold. Consequently $R$ is always a trivial line bundle and if one were to orient $R$ using $q$, the normal Hessian of $q$ along $D$ is positive definite. In other words, let $\alpha$ be a trivialization of $R^*$. Then $g:= \alpha(q) \in C^\infty(X)$ is a function with $g^{-1}(0) = D$ and ${\rm Hess}(g) = \alpha({\rm Hess}(q))$. Moreover, $D$ is a nondegenerate critical submanifold of $g$, and $g$ is locally Morse--Bott around $D$. As $D$ is codimension two, $X \setminus D$ is connected, so that the sign of $g$ on $X \setminus D$ is fixed. Replace $\alpha$ by $-\alpha$ if necessary so that this sign is positive, and then $q$ and $\alpha$ induce compatible orientations. Call such a trivialization $\alpha$ \emph{compatible} with $q$. For compatible trivializations we have $g \geq 0$ so that ${\rm Hess}(g)$ is positive definite. As a consequence of \autoref{lem:morsebott} we obtain the following.
\begin{prop}\label{prop:ellmorsebott} Let $|D| = (R,q)$ be an elliptic divisor and $\alpha$ a compatible trivialization. Then there exists a neighbourhood $U$ of the zero section $D \subset ND$ and an open embedding $\Phi\colon U \to X$ such that $\Phi|_D = {\rm id}_D$ and $(\Phi^*\alpha)(\Phi^* q) = Q_g$, where $g = \alpha(q) \in C^\infty(X)$.
\end{prop}
Given $p \in D$ we can locally trivialize the bundles $R$ and $ND$, so that using \autoref{prop:ellmorsebott} the section $q$ can be written locally as $q(x_1,\dots, x_n) = \pm(x_1^2 + x_2^2)$ in normal bundle coordinates such that $ND = \langle \partial_{x_1}, \partial_{x_2} \rangle$. Consequently, the elliptic ideal $I_D$ is locally generated by $r^2$, where $r^2 = x_1^2 + x_2^2$ is the squared radial distance from $D$ inside $ND$.
\begin{rem} While \autoref{prop:hypersurfacedivisor} shows that hypersurfaces carry a unique log divisor structure, the same is not true for codimension-two submanifolds and elliptic divisors. A simple example is provided by $X = \R^3$ with $D = \R \times \{0\} \times \{0\}$ and coordinates $(x,y,z)$. Equip $D$ with the elliptic ideals $I = \langle y^2 + z^2 \rangle$ and $I' = \langle y^2 + 2z^2 \rangle$. As these ideals are distinct, they supply $D$ with two non-isomorphic yet diffeomorphic elliptic divisor structures.
\end{rem}
\subsection{Divisors and geometric structures}
\label{sec:divgeomstr}
Divisors provide a convenient way to define and study specific classes of Poisson structures. We first describe the general way in which divisors relate to Poisson geometry. 

Let $X$ be a $2n$-dimensional manifold and $\pi$ a Poisson structure on $X$. The Pfaffian $\wedge^n \pi$ of $\pi$ is a section of the real line bundle $\wedge^{2n} TX$.
\begin{defn} A Poisson structure $\pi$ is of \emph{divisor-type} if $(\wedge^{2n} TX, \wedge^n \pi)$ is a divisor.
\end{defn}
Poisson structures of divisor-type are more thoroughly explored in \cite{Klaasse17}. In this paper we will restrict our attention to the two types associated to log and elliptic divisors.
\subsubsection{Log-Poisson and log-symplectic structures}
\label{sec:logsymp}
In this section we discuss Poisson structures built out of log divisors. We refer to \cite{CavalcantiKlaasse16, Klaasse17} for a more comprehensive treatment.
\begin{defn} Let $X$ be a $2n$-dimensional manifold. A \emph{log Poisson structure} is a Poisson structure $\pi$ on $X$ that is of log divisor-type.

\end{defn}
While the name ``log-Poisson structure'' is most consistent with other Poisson structures built out of divisors, these Poisson structures also go under the name of $b$-Poisson, $b$-symplectic, and \emph{log-symplectic structures} \cite{Cavalcanti13,GualtieriLi14,GuilleminMirandaPires14,MarcutOsornoTorres14,MarcutOsornoTorres14two}. The latter two names emphasize the ability to view log-Poisson structures as suitably degenerate symplectic forms, see \autoref{prop:blogbsymp}.

Let $(X,Z)$ be a \emph{log pair}, i.e.\ a manifold $X$ together with a log divisor $Z$. We say the pair $(X,Z)$ \emph{admits a log-Poisson structure} if there exists a log-Poisson structure $\pi$ on $X$ such that $Z_{\wedge^n \pi} = Z$. We will also denote $Z_{\wedge^n \pi}$ as $Z_\pi$.

The rank of a log-Poisson structure $\pi$ is equal to $2n$ on $X \setminus Z$, and $2n-2$ on $Z$. By the Weinstein splitting theorem \cite{Weinstein83}, any log-Poisson structure $\pi$ on $X^{2n}$ is locally equivalent around points in $Z_\pi$ to $x \partial_x \wedge \partial_y + \omega_0^{-1}$ on $\R^2 \times \R^{2n-2}$, where $\omega_0$ is the standard symplectic structure, and $Z_\pi = \{x = 0\}$ (see \cite{GuilleminMirandaPires14}). On $Z_\pi$, a \blog $\pi$ induces an equivalence class of \emph{cosymplectic structures}, which is a pair $(\alpha,\beta)$ of closed one- and two-forms on $Z$ such that $\alpha \wedge \beta^{n-1} \neq 0$. Here the equivalence relation is generated by addition of two-forms $df \wedge \alpha$, with $f \in C^\infty(Z)$.
Any compact cosymplectic manifold, such as $Z_\pi$, fibers over $S^1$ \cite{GualtieriLi14, Li08, MarcutOsornoTorres14, OsornoTorres15}. The following innocuous result will be important later in Section \ref{sec:boundarymaps}. Recall that a hypersurface $Z$ is \emph{separating} if $[Z] = 0 \in H_{2n-1}(X;\Z_2)$.
\begin{prop}\label{prop:orlogseparating} Let $(X,Z)$ be an orientable log-Poisson pair. Then $Z$ is separating.
\end{prop}
\bp As $X$ is orientable, $\wedge^{2n}TX$ is trivial. The section $\wedge^n \pi$ has transverse zeros, hence its zero set $Z$ is Poincar\'e dual to the first Stiefel--Whitney class of the trivial line bundle.
\ep
Any trivialization $\sigma \in \Gamma(\wedge^{2n}TX)$ produces a global defining function $h$ for $Z$ via $\wedge^n\pi = h \sigma$.
\subsubsection{Elliptic Poisson and stable generalized complex structures}
We can also construct Poisson structures out of elliptic divisors, obtaining the notion of an elliptic Poisson structure. These are called Poisson structures of elliptic log-symplectic type in \cite{CavalcantiGualtieri15}.
\begin{defn}[{\cite[Definition 3.3]{CavalcantiGualtieri15}}] Let $X$ be a $2n$-dimensional manifold. An \emph{elliptic Poisson structure} is a Poisson structure $\pi$ on $X$ that is of elliptic divisor-type.
\end{defn}
Let $(X,|D|)$ be an \emph{elliptic pair}, i.e.\ a manifold equipped with an elliptic divisor $|D|$. An elliptic pair $(X,|D|)$ \emph{admits an elliptic Poisson structure} if there exists an elliptic Poisson structure $\pi$ on $X$ such that $|D_{\wedge^n \pi}| = |D|$. When there is no elliptic divisor structure on $D$, we say that $(X,D)$ admits an elliptic Poisson structure if there exists some elliptic divisor structure $|D|$ on $D$ such that $(X,|D|)$ admits an elliptic Poisson structure.

The underlying Poisson structure of a \sgcs $\mc{J}$ can be shown to belong to this class, which in fact characterizes when $\mc{J}$ is stable.
\begin{prop}[{\cite[Theorem 3.7]{CavalcantiGualtieri15}}]\label{prop:stableellpoiss} Let $(X,H,\mc{J})$ be a generalized complex manifold. Then $\mc{J}$ is a \sgcs if and only if $\pi_{\mc{J}}$ is an elliptic Poisson structure.
\end{prop}
In Section \ref{sec:liealgbdsymp} we will see that this proposition can be pushed further. Namely, \sgcs{}s $\mc{J}$ on $X$ are in one-to-one correspondence with certain types of elliptic Poisson structures $\pi$ via the map $\mc{J} \mapsto \pi_{\mc{J}}$ (see \autoref{thm:sgcscorrespondence}). Moreover, the closed three-form $H$ required to state integrability of $\mc{J}$ is determined by $\pi$.

An elliptic pair $(X,|D|)$ \emph{admits a \sgcs} if there exists a closed three-form $H \in \Omega^3_{\rm cl}(X)$ and a \sgcs $\mc{J}$ on $(X,H)$ such that $|D_{\mc{J}}| = |D|$. When no elliptic divisor structure is specified on $D$, we say that $(X,D)$ admits a \sgcs if there exists some elliptic divisor structure $|D|$ on $D$ such that $(X,|D|)$ admits a \sgcs.
\newcommand{\ila}{ideal Lie algebroid\xspace}
\newcommand{\Ila}{Ideal Lie algebroid\xspace}
\section{Lie algebroids and Lie algebroid symplectic structures}
\label{sec:lalgebroids}
In this section we deal with some aspects of the general theory of symplectic forms in a Lie algebroid. For interesting applications, the Lie algebroid should be chosen so that such symplectic forms describe a certain type of geometric structure. The main examples relevant to this paper are the elliptic tangent bundle and the log-tangent bundle, describing \sgcs{}s and \blog{}s respectively. More examples of Lie algebroids with interesting symplectic geometry can be found in \cite{Klaasse17}.
\begin{defn} A \emph{Lie algebroid} is a vector bundle $\mc{A} \to X$ together with a Lie bracket $[\cdot,\cdot]_{\mc{A}}$ on $\Gamma(\mc{A})$ and a vector bundle map $\rho_{\mc{A}}\colon \mc{A} \to TX$ called the \emph{anchor}, satisfying the Leibniz rule $[v,f w]_{\mc{A}} = f [v,w]_{\mc{A}} + \mc{L}_{\rho_{\mc{A}}(v)} f \cdot w$ for all $v, w \in \Gamma(\mc{A})$, and $f \in C^\infty(X)$.
\end{defn}
We think of Lie algebroids as generalizations of the tangent bundle $TX$, chosen such that geometric constructions done using $\mc{A}$ are more suitable to the situation at hand. This should be kept in mind especially when we define Lie algebroid fibrations in Section \ref{sec:lalf}, as we will replace many usual notions by their Lie algebroid counterparts. The Lie algebroids we consider will be isomorphic via their anchor to the tangent bundle outside of a submanifold of positive codimension. We introduce the following subset to keep track of this.
\begin{defn} Let $\mc{A} \to X$ be a Lie algebroid. The \emph{isomorphism locus} of $\mc{A}$ is the open set $X_{\mc{A}} \subseteq X$ where $\rho_{\mc{A}}$ is an isomorphism.
\end{defn}
Any Lie algebroid gives rise to a graded algebra $\Omega^\bullet(\mc{A}) = \Gamma(\wedge^\bullet \mc{A}^*)$ of differential $\mc{A}$-forms. This algebra comes equipped with a differential $d_\mc{A}$ squaring to zero, constructed using the bracket $[\cdot,\cdot]_{\mc{A}}$ by means of the Koszul formula.
\begin{defn} Let $\mc{A} \to X$ be a Lie algebroid. The \emph{Lie algebroid cohomology} of $\mc{A}$ is given by $H^k_\mc{A}(X) = H^k(\Omega^\bullet(\mc{A}), d_{\mc{A}})$.
\end{defn}
For notational convenience we denote by $\Omega^k_\mc{A}(U)$, for an open subset $U \subset X$, the set of $\mc{A}$-$k$-forms defined on $U$, so that $\Omega^k_\mc{A}(X) = \Omega^k(\mc{A})$. The inclusion $i\colon X_{\mc{A}} \hookrightarrow X$ of the isomorphism locus induces a bijection $\rho_{\mc{A}}^*\colon \Omega^k(X_{\mc{A}}) \to \Omega^k_{\mc{A}}(X_{\mc{A}})$.
\subsection{The log-tangent bundle and elliptic tangent bundle}
\label{sec:logellbundles}
In this section we discuss the Lie algebroids that are of primary concern to us. These are the log-tangent bundle and the elliptic tangent bundle, constructed out of log divisors and elliptic divisors respectively. Both are examples of \emph{ideal Lie algebroids} \cite{Klaasse17}, which we now introduce.

Let $X$ be an $n$-dimensional manifold and let $I$ be an ideal sheaf, and denote by $\mc{V}_X$ the sheaf of vector fields on $X$. Let $\mc{V}(I) := \{v \in \mc{V}_X \, | \, \mc{L}_v I \subset I\}$ be the sheaf of derivations preserving $I$. This is a subsheaf of Lie algebras of $\mc{V}_X$. If $\mc{V}(I)$ is in addition locally finitely generated projective, it gives rise to a Lie algebroid $\mc{A}_I$ with $\Gamma(\mc{A}_I) = \mc{V}(I)$ by the Serre--Swan theorem.
\begin{defn} The Lie algebroid $\mc{A}_I$ is called the \emph{ideal Lie algebroid} associated to $I$.
\end{defn}
As $\mc{V}(I)$ is a submodule of $\Gamma(TX)$, the anchor of any ideal Lie algebroid is the natural inclusion on sections. However, as vector bundle map the anchor of $\mc{A}_I$ need not be an isomorphism. More precisely, the isomorphism locus $X_{\mc{A}_I}$ of $\mc{A}_I$ is the complement of the support ${\rm supp}(C^\infty(X) / I)$ of the quotient sheaf, or equivalently, the complement of ${\rm supp}(\mc{V}_X / \mc{V}(I))$. 
\subsubsection{The log-tangent bundle}
Let $Z = (L,s)$ be a log divisor on $X$ with associated ideal $I_Z = s (\Gamma(L^*))$. This ideal is exactly the vanishing ideal of $Z$, and $\mc{V}(I_Z)$ is the locally free sheaf of vector fields tangent to $Z$ \cite{Melrose93}. In local adapted coordinates $(z,x_2,\dots,x_n)$ around $Z = \{z = 0\}$ with $I_z = \langle z \rangle$, one has $\mc{V}(I_Z) = \langle z \partial_z, \partial_{x_2}, \dots, \partial_{x_n} \rangle$.

\begin{defn} The \emph{log-tangent bundle} $TX(- \log Z) \to X$ is the ideal Lie algebroid on $X$ with $\Gamma(TX(- \log Z)) = \mc{V}(I_Z)$.
\end{defn}
It is immediate that the isomorphism locus of $TX(-\log Z)$ is given by $X \setminus Z$. In analogy with the holomorphic case, denote $\Omega^k(\log Z) = \Omega^k(TX(-\log Z))$.
\begin{rem} The log-tangent bundle is also called the \emph{$b$-tangent bundle} \cite{Melrose93} and is then denoted by ${}^b TX$. We use the name log-tangent bundle as it shows the similarities with the elliptic tangent bundle defined below, and its notation allows us to keep track of $Z$.
\end{rem}
The log-tangent bundle admits a residue map. The restriction of $TX(-\log Z)$ to $Z$ surjects onto $TZ$ via the anchor map $\rho$, giving the following exact sequence
\be
	0 \to \mathbb{L}_Z \to TX(-\log Z)|_Z \stackrel{\rho}{\to} TZ \to 0,
\ee
where $\mathbb{L}_Z \to Z$ is the kernel of $\rho$. Dualizing gives a projection map ${\rm Res}_Z\colon \Omega^k(\log Z) \to \Omega^{k-1}(Z)$, which fits in the residue sequence
\be
	0 \to \Omega^\bullet(X) \stackrel{\rho^*}{\to} \Omega^\bullet(\log Z) \stackrel{{\rm Res}_Z}{\to} \Omega^{\bullet - 1}(Z) \to 0.
\ee
In terms of the local coordinate system above, with $\Gamma(TX(-\log Z)) = \langle z \partial_z, \partial_{x_2}, \dots, \partial_{x_n} \rangle$, a given log $k$-form $\alpha \in \Omega^k(\log Z)$ can be expressed as
\be
	\alpha = d \log z \wedge \alpha_0 + \alpha_1,
\ee
with $\alpha_i$ smooth forms. The inclusion $j_Z\colon Z \hookrightarrow X$ gives ${\rm Res}_Z(\alpha) = j_Z^* \alpha_0$. The following result referred to as the Mazzeo--Melrose theorem shows the above sequence splits, and identifies the Lie algebroid cohomology of the log-tangent bundle $TX(-\log Z)$ in terms of $X$ and $Z$.
\begin{thm}[\cite{Melrose93}]\label{thm:loglacohomology} Let $(X,Z)$ be a log pair. Then $H^k(TX(-\log Z)) \cong H^k(X) \oplus H^{k-1}(Z)$.
\end{thm}
\subsubsection{The elliptic tangent bundle}
Similarly, let $(X,|D|)$ be an elliptic pair with associated ideal $I_D = q (\Gamma(R^*))$. Note this is not the vanishing ideal of $D$. The associated submodule $\mc{V}(I_D)$ defines a sheaf of locally constant rank \cite{CavalcantiGualtieri15}, locally generated in appropriate polar coordinates $(r,\theta,x_3,\dots,x_n)$ around $D = \{r = 0\}$ such that $I_D = \langle r^2 \rangle$ by $\langle r \partial_r, \partial_\theta, \partial_{x_3}, \dots, \partial_{x_n} \rangle$.
\begin{defn} The \emph{elliptic tangent bundle} $TX(-\log|D|) \to X$ is the ideal Lie algebroid on $X$ with $\Gamma(TX(-\log|D|)) = \mc{V}(I_D)$.
\end{defn}
The isomorphism locus of $TX(-\log |D|)$ is given by $X \setminus D$. As for the log-tangent bundle, we denote $\Omega^k(\log |D|) = \Omega^k(TX(-\log |D|))$.
\begin{rem} There is another ideal Lie algebroid defined in \cite{CavalcantiGualtieri15}, namely the \emph{complex log-tangent bundle} $TX(-\log D)$ (there called the \emph{logarithmic tangent bundle}). These are associated to complex log divisors $D = (U,\sigma)$. This Lie algebroid is a generalization of the log-tangent bundle one can define on complex manifolds equipped with a divisor. We will not directly use this Lie algebroid in this paper, hence we will not elaborate on its properties.
\end{rem}

The elliptic tangent bundle admits several residue maps \cite{CavalcantiGualtieri15}, three of which we will now describe. The \emph{elliptic residue} ${\rm Res}_q$ comes from considering the restriction of $TX(-\log |D|)$ to $D$, which fits in an exact sequence
\begin{equation}
\label{eqn:ellipticrestriction}
	0 \to \underline{\R} \oplus \mf{k} \to TX(-\log |D|)|_D \to TD \to 0,
\end{equation}
with $\underline{\R}$ generated by the Euler vector field of $ND$, and $\mf{k} \cong \wedge^2 N^*D \otimes R$. Choosing a coorientation for $D$, i.e.\ a trivialization of $ND$, also trivializes $\mf{k}$. Dualizing the above sequence we obtain a projection map ${\rm Res}_q\colon \Omega^k(\log |D|) \to \Omega^{k-2}(D,\mf{k}^*)$, with kernel $\Omega^\bullet_0(\log |D|)$ the subcomplex of $\Omega^\bullet(\log |D|)$ of zero elliptic residue forms. 

Denote by $S^1ND$ the circle bundle associated to $ND$. The \emph{radial residue} ${\rm Res}_r$ arises from quotienting \eqref{eqn:ellipticrestriction} by the Euler vector field of $ND$, giving the extension
\be
	0 \to \mf{k} \to {\rm At}(S^1 ND) \to TD \to 0,
\ee
where ${\rm At}(S^1 ND)$ is the associated Atiyah algebroid of $S^1 ND$. Noting that $TX(-\log |D|)|_D$ is a trivial extension of ${\rm At}(S^1 ND)$, the elliptic residue factors through the radial residue map ${\rm Res}_r\colon \Omega^k(\log |D|) \to \Gamma(D, \wedge^{k-1} {\rm At}(S^1  ND)^*)$. When the elliptic residue vanishes, the radial residue naturally maps to $\Omega^{k-1}(D)$ without needing a coorientation. Finally, there is a \emph{$\theta$-residue} ${\rm Res}_\theta:\, \Omega^k_0(\log|D|) \to \Omega^{k-1}(D)$, which we will only define for forms with zero elliptic residue. We provide a description of these residue maps in local coordinates. In the adapted coordinate system around $D$ as above, where $\Gamma(TX(-\log |D|)) = \langle r \partial_r, \partial_\theta, \partial_{x_3}, \dots, \partial_{x_n} \rangle$, a given elliptic $k$-form $\alpha \in \Omega^k(\log |D|)$ can be locally written as
\be
	\alpha = d \log r \wedge d\theta \wedge \alpha_0 + d \log r \wedge \alpha_1 + d\theta \wedge \alpha_2 + \alpha_3,
\ee
where each $\alpha_i$ is a smooth form. Using the inclusion $j_D\colon D \hookrightarrow X$ we have ${\rm Res}_q(\alpha) = j_D^* \alpha_0$, ${\rm Res}_r(\alpha) = (d\theta \wedge \alpha_0 + \alpha_1)|_D$, and we set ${\rm Res}_\theta(\alpha) = j_D^*\alpha_2$. Moreover, we see that ${\rm Res}_q(\alpha) = \iota_{\partial_\theta} {\rm Res}_r(\alpha)$.

As for the log-tangent bundle, the Lie algebroid cohomology of the elliptic tangent bundle as well as its zero elliptic residue version can be expressed in terms of $X$ and $D$.
\begin{thm}[{\cite[Theorems 1.8, 1.12]{CavalcantiGualtieri15}}]\label{thm:elllacohomology} Let $(X,|D|)$ be an elliptic pair. Then $H^k(\log |D|) \cong H^k(X \setminus D) \oplus H^{k-1}(S^1 ND)$. Further, one has $H^k_0(\log |D|) \cong H^k(X \setminus D) \oplus H^{k-1}(D)$.
\end{thm}
The above isomorphisms are induced by the maps $(i^*, {\rm Res}_r)$ with $i\colon X\setminus D \hookrightarrow X$ the inclusion of the divisor complement, noting ${\rm Res}_r$ naturally maps to $\Omega^{k-1}(D)$ when the elliptic residue vanishes (in that case, ${\rm Res}_r(\alpha) = \alpha_1|_D$ in the above local coordinates).
\subsection{Lie algebroid morphisms}
In this section we discuss morphisms between Lie algebroids, focusing specifically on those between the elliptic and log-tangent bundle. 
\begin{defn} Let $\mc{A}, \mc{A}' \to X$ be Lie algebroids over the same base $X$. A \emph{Lie algebroid morphism} is a vector bundle map $(\varphi, f)\colon \mc{A} \to \mc{A}'$ such that $df \circ \rho_{\mc{A}} = \rho_{\mc{A}'} \circ \varphi$ and $\varphi[v,w]_{\mc{A}} = [\varphi(v),\varphi(w)]_{\mc{A}'}$ for all $v,w \in \Gamma(\mc{A})$.
\end{defn}
The above definition does not immediately generalize to varying base, as a vector bundle map does not give a map on the space of sections, as is required in order to state the bracket condition. There is a description of bracket compatibility in terms of the pullback bundle $f^* \mc{A}'$, but we will instead use the following equivalent definition in terms of their duals.

Let $\mc{A} \to X$ and $\mc{A}' \to X'$ be Lie algebroids. Vector bundle maps $(\varphi,f)\colon \mc{A} \to \mc{A}'$ are in one-to-one correspondence with algebra morphisms $\varphi^*\colon \Omega^\bullet(\mc{A}') \to \Omega^\bullet(\mc{A})$. Using this we can phrase the conditions that $\varphi$ preserves anchors and brackets in terms of $\varphi^*$.
\begin{defn} A vector bundle map $(\varphi,f)\colon \mc{A} \to \mc{A}'$ is a \emph{Lie algebroid morphism} if $\varphi^*\colon \Omega^\bullet(\mc{A}') \to \Omega^\bullet(\mc{A})$ is a chain map, that is $d_{\mc{A}'} \circ \varphi^* = \varphi^* \circ d_\mc{A}$.
\end{defn}
A Lie algebroid morphism $(\varphi,f)\colon \mc{A} \to \mc{A}'$ is in particular a morphism of anchored vector bundles, i.e.\ satisfies $df \circ \rho_\mc{A} = \rho_{\mc{A}'} \circ \varphi$. Moreover, a Lie algebroid morphism restricts pointwise to a linear map $\varphi\colon \ker \rho_{\mc{A}} \to \ker \rho_{\mc{A}'}$ between kernels of the respective anchor maps. Furthermore, $\rho_\mc{A}$ gives an isomorphism between $\ker \varphi$ and $\ker df$ when in the isomorphism locus.
For the Lie algebroids we will consider, smooth maps of the underlying manifolds give rise to Lie algebroid morphisms, as long as they intertwine the anchor maps. This is true in general for anchored vector bundle morphisms between Lie algebroids with dense isomorphism loci.
\begin{prop}\label{prop:denselam} Let $\mc{A}_X \to X$ and $\mc{A}_Y \to Y$ be Lie algebroids such that $X_{\mc{A}_X}$ is dense. Suppose that $(\varphi,f)\colon \mc{A}_X \to \mc{A}_Y$ is an anchored bundle morphism and $f^{-1}(Y_{\mc{A}_Y}) = X_{\mc{A}_X}$. Then $(\varphi, f)\colon \mc{A}_X \to \mc{A}_Y$ is a Lie algebroid morphism.
\end{prop}
\bp As $\varphi$ is a vector bundle morphism, $\varphi^*$ is an algebra morphism. In the isomorphism loci, $\varphi$ must equal $df$, and $(df,f)$ is a Lie algebroid morphism between $TX$ and $TY$, i.e.\ $f^*$ is a chain map. By density of $X_{\mc{A}_X}$, the map $\varphi^*$ is a chain map everywhere.
\ep
Consequently, for such Lie algebroids one can determine whether there is a Lie algebroid morphism $(\varphi, f)$ by checking that $f^*$ extends to a map $\varphi^*$ on forms. This in turn will follow by the universal property of the exterior algebra if it holds on generators, so it suffices to check that $f^*$ extends to $\varphi^*\colon \Omega^1(\mc{A}_Y) \to \Omega^1(\mc{A}_X)$. We will now perform such a check to obtain a Lie algebroid morphism between the elliptic and log tangent bundles.
\subsubsection{Morphisms between elliptic and log-tangent bundles}
Let $(X, |D|)$ be an elliptic pair and $(Y,Z)$ a log pair, with associated ideals $I_D \subset C^\infty(X)$ and $I_Z \subset C^\infty(Y)$. A morphism of the corresponding divisors gives rise to a Lie algebroid morphism.
\begin{prop}\label{prop:elltologmorphism} Let $f\colon (X,|D|) \to (Y,Z)$ be a morphism of divisors. Then $df$ induces a Lie algebroid morphism $(\varphi,f)\colon TX(-\log|D|) \to TY(-\log Z)$ such that $\varphi = df$ on sections.
\end{prop}
\bp It is immediate that $f^{-1}(Z) = D$. Note that the isomorphism loci of the elliptic tangent bundle and the log-tangent bundle are dense. By \autoref{prop:denselam} it thus suffices to show that $df$ induces a vector bundle morphism, which in turn is equivalent to showing that $f^*$ extends to an algebra morphism from $\Omega^\bullet(Y;\log Z)$ to $\Omega^\bullet(X; \log |D|)$. Let $x \in D$ and denote $y = f(x) \in Z$. Consider suitable tubular neighbourhood coordinates $(r, \theta, x_3,\dots, x_n)$ in a neighbourhood $U$ of $x$ such that $U \cap D = \{r = 0\}$ with $I_D = \langle r^2 \rangle$, and $(z, y_2, \dots, y_m)$ in a neighbourhood $V$ of $y$ such that $V \cap Z = \{z = 0\}$ and $I_Z = \langle z \rangle$. In these coordinates we have $\Omega^1(U; - \log |D|) = \langle d \log r, d \theta, dx_3, \dots, d x_n\rangle$ and $\Omega^1(V; -\log Z) = \langle d\log z, dy_2, \dots, dy_m\rangle$.

The Lie algebroid one-forms $dy_i$ can be pulled back using $f^*$ as these are smooth. Moreover, the smooth one-forms inject into $\Omega^1(U;- \log D)$ using the anchor. We are left with checking that $d \log z$ is pulled back to a Lie algebroid one-form. As $f$ is a morphism of divisors we have $f^* I_Z = I_D$, so that $f^*(z) = e^h r^2$ for $h$ a smooth function $U$. Consequently $f^* d \log z = d \log f^*(z) = d \log (e^h r^2) = d h + 2 d\log r \in \Omega^1(U;-\log |D|)$. We conclude that $df$ induces a Lie algebroid morphism as desired.
\ep
\begin{rem}\label{rem:otherlamorphisms} Lie algebroid morphisms between log-tangent bundles arise out of so-called \emph{$b$-maps}, which are maps $f\colon (X,Z_X) \to (Y,Z_Y)$ between log pairs such that $f^{-1}(Z_Y) = Z_X$ and $f$ is transverse to $Z_Y$. Alternatively, they are exactly the corresponding morphisms of divisors, i.e.\ those maps satisfying $f^* I_{Z_Y} = I_{Z_X}$ \cite{Klaasse17}. Similarly, morphisms between elliptic tangent bundles arise from those maps $f\colon (X,|D_X|) \to (Y,|D_Y|)$ between elliptic pairs for which $f^{-1}(D_Y) = D_X$ and $f$ is transverse to $D_Y$. Again, these can be alternatively described as maps satisfying $f^* I_{D_Y} = I_{D_X}$ \cite{Klaasse17}.
\end{rem}
\subsubsection{Lie algebroid morphisms and residue maps}
A Lie algebroid morphism from the elliptic tangent bundle $TX(-\log |D|)$ to the log tangent bundle $TY(-\log Z)$ intertwines the residue maps that were discussed in Section \ref{sec:logellbundles}.
\begin{prop}\label{prop:residuemaps} Let $(\varphi,f)\colon TX(-\log |D|) \to TY(-\log Z)$ be a Lie algebroid morphism. Then ${\rm Res}_q \circ \varphi^* = 0$. Moreover, $({\rm Res}_r+{\rm Res}_\theta) \circ \varphi^* = f^* \circ {\rm Res}_Z$.
\end{prop}
To prove the above proposition, we briefly elaborate on how the residue maps come about. Given a short exact sequence $\mc{S}\colon 0 \to E \to W \to V \to 0$ of vector spaces, there is an associated dual sequence $\mc{S}^*\colon 0 \to V^* \to W^* \to E^* \to 0$. For a given $k \in \N$, by taking $k$th exterior powers we obtain a filtration of spaces $\mc{F}^i := \{\rho \in \wedge^k W^* \, | \, \iota_x \rho = 0$ for all $x \in \wedge^i E \}$, for $i = 0,\dots,k+1$. These spaces satisfy $\mc{F}^0 = 0$, $\mc{F}^1 = \wedge^k V^*$, $\mc{F}^i \subset \mc{F}^{i+1}$, and $\mc{F}^{i+1} / \mc{F}^i \cong \wedge^{k-i} V^* \otimes \wedge^i E^*$. Setting $\ell := \dim E$, we have $\mc{F}^{\ell + 1} = \wedge^k W^*$.
\begin{defn} The \emph{residue} of $\rho \in \wedge^k W^*$ is ${\rm Res}(\rho) = [\rho] \in \mc{F}^{\ell + 1} / \mc{F}^{\ell} \cong \wedge^{k-\ell} V^* \otimes \wedge^\ell E^*$.
\end{defn}
Upon a choice of trivialization of $\wedge^\ell E^*$, i.e.\ a choice of orientation for $E$, one can view the residue ${\rm Res}(\rho)$ as an element of $\wedge^{k-\ell} V^*$. One can further consider lower residues ${\rm Res}_{-m}\colon \wedge^k W^* \to \mc{F}^{\ell+1} / \mc{F}^{\ell-m}$ for $m > 0$. These are always defined but have a better description for forms $\rho \in \wedge^k W^*$ whose higher residues vanish, so that ${\rm Res}_{-m}(\rho) \in \mc{F}^{\ell-m +1} / \mc{F}^{\ell -m} \cong \wedge^{k-\ell+m} V^* \otimes \wedge^{\ell-m} E^*$. Given a map of short exact sequences $\Psi\colon \wt{\mc{S}} \to \mc{S}$ with dual map $\Psi^*\colon \mc{S}^* \to \wt{\mc{S}}^*$, there is a corresponding map of filtrations $\Psi^*\colon \mc{F}^i \to \wt{\mc{F}}^i$. Setting $\wt{\ell} := \dim \wt{E}$, we have the following.
\begin{lem}\label{lem:residues} In the above setting, assume that $\wt{\ell} > \ell$. Then $\wt{\rm Res}(\Psi^* \rho) = 0$ for all $\rho \in \wedge^k W^*$.
\end{lem}
\bp We have $\rho \in \mc{F}^{\ell +1}$ so that $\Psi^* \rho \in \wt{\mc{F}}^{\ell + 1}$. As $\wt{\ell} > \ell$, we have $\wt{\mc{F}}^{\ell+1} \subset \wt{\mc{F}}^{\wt{\ell}} \subset \wt{\mc{F}}^{\wt{\ell}+1}$, so that $\wt{\rm Res}(\Psi^* \rho) = [\Psi^* \rho] \in \wt{\mc{F}}^{\wt{\ell}+1} / \wt{\mc{F}}^{\wt{\ell}}$ vanishes as desired.
\ep
Assuming $\wt{\ell} > \ell$, all lower residues automatically vanish by degree reasons until considering $\Psi^* \rho \in \wt{\mc{F}}^{\ell +1}$. Hence the first possibly nonzero residue is $\wt{\rm Res}_{\ell-\wt{\ell}}(\Psi^* \rho) = [\Psi^* \rho] \in \wt{\mc{F}}^{\ell+1} / \wt{\mc{F}}^{\ell} \cong \wedge^{k-\ell} \wt{V}^* \otimes \wedge^\ell \wt{E}^*$. Similarly to \autoref{lem:residues} we obtain the following.
\begin{lem}\label{lem:residuecommute} In the above setting, assuming $\wt{\ell} \geq \ell$, we have $\Psi^* \circ {\rm Res} = \wt{\rm Res}_{\ell - \wt{\ell}} \circ \Psi^*$.
\end{lem}
With this understood, we can turn to proving \autoref{prop:residuemaps}.
\bp[ of \autoref{prop:residuemaps}] By assumption we have $f^{-1}(Z) = D$ so that $df\colon TD \to TZ$. Restricting $TX(-\log |D|)$ to $D$ and $TY(-\log Z)$ to $Z$ gives the following commutative diagram.
\begin{center}
	\begin{tikzpicture}

	\matrix (m) [matrix of math nodes, row sep=2.5em, column sep=2.5em,text height=1.5ex, text depth=0.25ex]
	{	\wt{\mc{S}}\colon 0 & \underline{\R} \oplus \mf{k} & TX(-\log |D|)|_D & TD & 0 \\ \mc{S}\colon 0 & \mathbb{L}_Z & TY(-\log Z)|_Z & TZ & 0 \\};
	\path[-stealth]
	(m-1-1) edge (m-1-2)
	(m-1-2) edge (m-1-3)
	(m-1-4) edge (m-1-5)
	(m-2-1) edge (m-2-2)
	(m-2-2) edge (m-2-3)
	(m-2-4) edge (m-2-5)
	(m-1-3) edge node [left] {$\varphi$} (m-2-3)
	(m-1-2) edge node [left] {$\varphi$} (m-2-2)
	(m-1-4) edge node [right] {$df$} (m-2-4)
	(m-1-3) edge node [above] {$\rho_X$} (m-1-4)
	(m-2-3) edge node [above] {$\rho_Y$} (m-2-4);	
	\end{tikzpicture}
\end{center}
Consequently, we obtain a map $\varphi^*\colon \mc{S}^* \to \wt{\mc{S}}^*$ between dual sequences, and also between spaces of sections. Using the notation preceding this proof we have $E = \mathbb{L}_Z$ so that $\ell = \dim(E) = 1$, and $\wt{E} = \underline{\R} \oplus \mf{k}$ so that $\wt{\ell} = \dim(\wt{E}) = 2$. Recall that $\mathbb{L}_Z$ carries a canonical trivialization. Given a form $\alpha \in \Omega^k(\log Z)$, we can identify ${\rm Res}(\alpha) \in \Gamma(Z;\wedge^{k-1} T^*Z \otimes \mathbb{L}_Z^*)$ with ${\rm Res}_Z(\alpha) \in \Omega^{k-1}(Z)$. Similarly, a choice of coorientation for $ND$ trivializes $\mf{k} = \wedge^2 N^*D \otimes R$. Given $\beta \in \Omega^k(\log |D|)$ this identifies $\wt{\rm Res}(\beta) \in \Gamma(D; \wedge^{k-2} T^*D \otimes \mf{k}^*)$ with ${\rm Res}_q(\beta) \in \Omega^{k-2}(D)$, using that $\wedge^2(\underline{\mathbb{R}} \oplus \mf{k}) \cong \mf{k}$. Moreover, for $\beta$ with ${\rm Res}_q(\beta) = 0$, the radial residue ${\rm Res}_r(\beta) \in \Omega^{k-1}(D)$ together with the $\theta$-residue ${\rm Res}_\theta(\beta) \in \Omega^{k-1}(D)$ is identified with $\wt{\rm Res}_{-1}(\beta) \in \Omega^{k-1}(D;\underline{\R} \oplus \mf{t})$. As $2 = \wt{\ell} \geq \ell = 1$ with $\ell - \wt{\ell} = -1$ we obtain immediately from \autoref{lem:residues} and \autoref{lem:residuecommute} that ${\rm Res}_q(\varphi^* \alpha) = 0$ and that $({\rm Res}_r+{\rm Res}_\theta)(\varphi^* \alpha) = f^*({\rm Res}_Z(\alpha))$.
\ep
\begin{rem} The proof of \autoref{prop:residuemaps} shows that the residue maps are also intertwined for morphisms between log-tangent bundles, and between elliptic tangent bundles.
\end{rem}
\subsection{Lie algebroid Lefschetz fibrations}
\label{sec:lalf}
In this section we introduce the appropriate notions of fibration and Lefschetz fibration in the context of Lie algebroid morphisms.
\begin{defn} A \emph{Lie algebroid fibration} is a Lie algebroid morphism $(\varphi,f)\colon \mc{A} \to \mc{A}'$ for which $\varphi\colon \mc{A} \to f^* \mc{A}'$ is surjective. Equivalently, $\varphi$ should be fiberwise surjective.
\end{defn}
Note that if $f\colon X \to Y$ is a fibration, then $(df,f)\colon TX \to TY$ is a Lie algebroid fibration.
%
\begin{rem} Our notion of a Lie algebroid fibration differs from the one used by other authors, notably Mackenzie \cite{Mackenzie05}. Our Lie algebroid fibrations are not required to cover a surjective submersion. In other words, only $\varphi$ should be fiberwise surjective, not both $\varphi$ and $df$. This is in line with viewing $\mc{A}$ as the replacement of $TX$.
\end{rem}
We next introduce Lie algebroid \lf{}s, which are a simultaneous generalization of \lf{}s, as well as of Lie algebroid fibrations.  We first recall the notion of a \lf. See \cite{GompfStipsicz99} for more information.
\begin{defn}\label{defn:lefschetzfibration} A \emph{\lf} is a proper map $f\colon X^{2n} \to \Sigma^2$ between oriented connected manifolds which is injective on critical points and such that for each critical point $x \in X$ there exist orientation preserving complex coordinate charts centered at $x$ and $f(x)$ in which $f$ takes the form $f\colon \C^n \to \C$, $f(z_1, \dots, z_n) = z_1^2 + \dots + z_n^2$.
\end{defn}
\begin{rem} The requirement that \lf{}s are injective on critical points is convenient but not essential. If it is not satisfied it can be ensured by a small perturbation.
\end{rem}
Next we essentially separate two types of singular behavior, namely that of the anchors of the Lie algebroids, and that of the morphism between them. In the isomorphism locus of the Lie algebroid, the condition of being a Lie algebroid fibration is just that it be a fibration. We can weaken this condition here, and only here, to allow for Lefschetz-type singularities.
\begin{defn} A \emph{Lie algebroid \lf} $(\varphi,f)\colon \mc{A}^{2n}_{X} \to \mc{A}^2_{\Sigma}$ is a Lie algebroid morphism which is a Lie algebroid fibration outside a discrete set $\Delta \subset X_{{\mc{A}}_X}$ with $f(\Delta) \subset \Sigma_{\mc{A}_\Sigma}$ such that $f|_{X_{\mc{A}_X}}\colon X_{{\mc{A}}_X} \to \Sigma_{{\mc{A}}_\Sigma}$ is a \lf.
\end{defn}
\begin{rem} Note that if $X_{\mc{A}_X}$ is empty, the notions of Lie algebroid \lf and fibration coincide. There are no Lefschetz singularities outside of $X_{\mc{A}_X}$. Moreover, whenever $\Delta$ is nonempty, both $X_{\mc{A}_X}$ and $\Sigma_{\mc{A}_\Sigma}$ are nonempty, and hence $\dim(X) = 2n$ and $\dim(\Sigma) = 2$.
\end{rem}
In Section \ref{sec:gompfthurston}, we will use Lie algebroid \lf{}s whose generic fibers in $X_{\mc{A}_X}$ are connected. Unlike for usual Lefschetz fibrations \cite{GompfStipsicz99}, for Lie algebroid Lefschetz fibrations there is in general no exact sequence in homotopy by which we can ensure connected fibers.
\subsection{Lie algebroid symplectic structures}
\label{sec:liealgbdsymp}
In this section we discuss symplectic structures for a Lie algebroid. We start by defining Poisson structures for Lie algebroids, as the symplectic forms we consider in this paper will be obtained through their bivector counterparts. Throughout, let $\mc{A} \to X$ be a Lie algebroid. The bracket $[\cdot,\cdot]_\mc{A}$ extends in the natural way to an $\mathcal{A}$-Schouten bracket on $\Gamma(\wedge^\bullet \mathcal{A})$, again denoted by $[\cdot,\cdot]_\mc{A}$.
\begin{defn} An \emph{$\mathcal{A}$-Poisson structure} is an $\mathcal{A}$-bivector $\pi_\mathcal{A}\in \Gamma(\wedge^2\mathcal{A})$ with $[\pi_\mathcal{A},\pi_\mathcal{A}]_{\mc{A}} = 0$.
\end{defn}
\begin{defn} A Poisson structure $\pi$ on $X$ is said to be \emph{of $\mc{A}$-type} for a Lie algebroid $\rho_\mc{A}\colon \mc{A} \to X$ if there exists an $\mc{A}$-bivector $\pi_{\mc{A}} \in \Gamma(\wedge^2 \mc{A})$ such that $\rho_{\mc{A}}(\pi_{\mc{A}}) = \pi$.
\end{defn}
In the above situation we call $\pi_{\mc{A}}$ an \emph{$\mc{A}$-lift} of $\pi$. Such a lift need not be unique if it exists (consider the zero Poisson structure and any $\mc{A}$ with zero anchor), but it is unique if $X_\mc{A}$ is dense, as will be the case for us. Moreover, when $X_{\mc{A}}$ is dense, any $\mc{A}$-lift of a Poisson structure is automatically $\mc{A}$-Poisson. Note that we have $\pi^\sharp = \rho_{\mc{A}} \circ \pi_{\mc{A}}^\sharp \circ \rho_{\mc{A}}^*$ as maps.
\begin{defn} An \emph{$\mc{A}$-symplectic structure} is a closed and nondegenerate $\mc{A}$-two-form ${\omega_{\mc{A}} \in \Omega^2_{\mc{A}}(X)}$. Denote the space of $\mc{A}$-symplectic forms by ${\rm Symp}(\mc{A})$.
\end{defn}
An $\mc{A}$-Poisson structure $\pi_{\mc{A}}$ is called \emph{nondegenerate} if $\pi_{\mc{A}}^\sharp$ is an isomorphism. As when $\mc{A} = TX$, for any Lie algebroid of even rank there is a bijection between $\mc{A}$-symplectic forms and nondegenerate $\mc{A}$-Poisson structures. Namely, given an $\mc{A}$-symplectic structure $\omega_{\mc{A}}$, nondegeneracy implies we can invert the map $\omega_{\mc{A}}^\flat\colon \mc{A} \to \mc{A}^*$ to $(\omega_{\mc{A}}^\flat)^{-1} = \pi_{\mc{A}}^\sharp\colon \mc{A}^* \to \mc{A}$ for an $\mc{A}$-Poisson structure $\pi_{\mc{A}}$. The conditions $d_{\mc{A}} \omega_{\mc{A}} = 0$ and $[\pi_{\mc{A}}, \pi_{\mc{A}}]_{\mc{A}} = 0$ are equivalent.

Given an $\mc{A}$-\symp{} $\omega_\mc{A}$, we call $\pi = \rho_{\mc{A}}(\pi_\mc{A})$ the \emph{dual bivector} to $\omega_{\mc{A}}$. We say $\pi$ is \emph{of nondegenerate $\mc{A}$-type} if it admits a nondegenerate $\mc{A}$-lift $\pi_{\mc{A}}$. In this paper we will focus on nondegenerate $\mc{A}$-Poisson structures, as we wish to use symplectic techniques.

The elliptic Poisson structures from Section \ref{sec:divgeomstr} are in one-to-one correspondence with symplectic forms in the associated elliptic tangent bundle, i.e.\ \emph{elliptic symplectic structures}.
\begin{prop}[{\cite[Lemma 3.4]{CavalcantiGualtieri15}}]\label{prop:elllogcorrespondence} A Poisson structure $\pi$ on $X^{2n}$ is elliptic if and only if it is of nondegenerate $\mc{A}$-type, where $\mc{A} = TX(-\log |D|)$ with elliptic divisor $|D| = (\wedge^{2n} TX, \wedge^n \pi)$.
\end{prop}
We can now state the extension of \autoref{prop:stableellpoiss}, giving a characterization of \sgcs{} purely in terms of elliptic symplectic structures.
\begin{thm}[{\cite[Theorem 3.7]{CavalcantiGualtieri15}}]\label{thm:sgcscorrespondence} Let $X$ be a compact manifold. There is a bijection $(\mc{J}, H) \to (\pi_{\mc{J}}^{-1}, \mf{o})$ between \sgcs{}s up to gauge equivalence and elliptic symplectic structures with vanishing elliptic residue and cooriented degeneracy locus.
\end{thm}
The associated closed three-form $H$ in the definition of a \gcs can be determined via $[H] = {\rm Res}_r([\omega_{\mc{A}}]) \wedge {\rm PD}[D]$, where $\omega_{\mc{A}}$ is the Lie algebroid symplectic structure for $\mc{A} = TX(-\log |D_{\mc{J}}|)$ whose dual bivector is $\pi_{\mc{J}}$, and $D = (\wedge^n \pi_\mc{J})^{-1}(0)$. The Poincar\'e dual of $D$ requires a choice of coorientation.

Similarly to \autoref{prop:elllogcorrespondence}, there is a bijective correspondence between log-Poisson structures, and \blog{}s. This was used in \cite{CavalcantiKlaasse16} to construct log-Poisson structures.
\begin{prop}[{\cite[Proposition 20]{GuilleminMirandaPires14}}]\label{prop:blogbsymp} A Poisson structure $\pi$ on $X^{2n}$ is log if and only if it is of nondegenerate $\mc{A}$-type, where $\mc{A} = TX(- \log Z)$ with log divisor $Z = (\wedge^{2n} TX, \wedge^n \pi)$.
\end{prop}
Finally, we note that a Lie algebroid $\mc{A} \to X$ of rank two admits an $\mc{A}$-\symp if and only if $\mc{A}$ is orientable, i.e.\ $w_1(\mc{A}) = 0$. Namely, any nonvanishing section of $\wedge^2 \mc{A}$ is a nondegenerate $\mc{A}$-Poisson structure by dimension reasons. For later use we record the following, also noted in \cite{CavalcantiKlaasse16}.
\begin{lem}\label{lem:blogsurface} Let $\Sigma^2$ be a compact oriented surface. Then $(\Sigma,\partial \Sigma)$ admits a \blog. For a hypersurface $Z$, $(\Sigma,Z)$ admits a \blog if and only if $[Z] = 0 \in H_1(\Sigma;\Z_2)$.
\end{lem}
\bp By the previous discussion, the pair $(\Sigma, Z)$ carries a \blog if and only if $T\Sigma(-\log Z)$ is orientable. Note that $w_1(T\Sigma(-\log Z)) = w_1(\Sigma) + w_1(L)$ via the bundle isomorphism $T\Sigma(-\log Z) \oplus L \cong T\Sigma \oplus \R$ \cite{Klaasse17}, where $L = L_{I_Z}$ in the notation of \autoref{prop:locprincideal}. As an orientable manifold with boundary has orientable boundary, the result follows, noting that $w_1(L) = {\rm PD}_{\Z_2}(Z)$, the Poincar\'e dual of $[Z]$ with $\Z_2$-coefficients.
\ep
\section{Constructing Lie algebroid symplectic structures}
\label{sec:gompfthurston}
In this section we consider the Thurston argument \cite{Thurston76} for constructing \symp{}s extended by Gompf \cite{Gompf04}, and adapt it to the context of Lie algebroid symplectic forms. The case of the log-tangent bundle can be found in \cite{CavalcantiKlaasse16}. The guiding principle is to combine suitably symplectic-type structures from the base of a fibration-like map with a form that is symplectic on the tangent spaces of the fibers of that map. In the Lie algebroid case one uses Lie algebroid morphisms $(\varphi,f)\colon \mc{A}_X \to \mc{A}_Y$. Special attention is required because $\rho_{\mc{A}_X}\colon \ker \varphi \to \ker df$ need not be an isomorphism (nor injective or surjective), hence one should interpret the tangent space to the fibers suitably.

We will use Lie algebroid \acs{}s as certificates for nondegeneracy of forms, by using the notion of tameness. Let $\mc{A}_X \to X$ be a Lie algebroid and $\omega \in {\rm Symp}(\mc{A}_X)$.
\begin{defn}\label{defn:taming} An \emph{$\mc{A}_X$-\acs} is a vector bundle complex structure $J$ for $\mc{A}_X$. An $\mc{A}_X$-\acs $J$ is \emph{$\omega$-tame} if $\omega(v, J v) > 0$ for all $v \in \Gamma(\mc{A}_X)$. Given $(\varphi,f)\colon \mc{A}_X \to \mc{A}_Y$ a Lie algebroid morphism and $\omega_Y \in {\rm Symp}(\mc{A}_Y)$, $J$ is \emph{$(\omega_Y,\varphi)$-tame} if $(\varphi_x ^* \omega_Y)(v, J v) > 0$ for all $v \in \mc{A}_{X,x} \setminus \ker \varphi_x$ and $x \in X$.
\end{defn}
As usual, the space of taming $\mc{A}_X$-\acs{}s for $\omega$ is convex and nonempty and is denoted by $\mc{J}(\omega)$. Note that any $\omega \in \Omega_{\mc{A}_X}^2(X)$ taming a $\mc{A}_X$-\acs $J$ is necessarily nondegenerate. Hence if $\omega$ is a closed $\mc{A}_X$-two-form taming some $J$, then $\omega$ is $\mc{A}_X$-symplectic and $J$ induces the same $\mc{A}_X$-orientation as $\omega$. Note moreover that if $J$ is $(\omega_Y, \varphi)$-tame, then $\ker \varphi$ is a $J$-complex subspace of $\mc{A}_X$. Indeed, if $v \in \ker \varphi$ and $J v \not\in \ker \varphi$, we would have $0 = \varphi^*\omega_Y(v, Jv) = \varphi^*\omega_Y(Jv, J(Jv)) > 0$, which is a contradiction.
\begin{prop}\label{prop:tamingopen} The taming condition is open, i.e.\ it is preserved under sufficiently small perturbations of $\omega$ and $J$, and of varying the point in $X$.
\end{prop}
\bp The taming condition $\omega(v, J v) > 0$ for the pair $(\omega, J)$ holds provided it holds for all $v \in \Sigma\mc{A}_X \subset \mc{A}_X$, the unit sphere bundle with respect to some preassigned metric, as $X$ is compact. As $\Sigma\mc{A}_X$ is also compact, the continuous map $\wt{\omega}\colon \Sigma\mc{A}_X \to \R$ given by $\wt{\omega}(v) := \omega(v, J v)$ for $v \in \Sigma\mc{A}_X$ is bounded from below by a positive constant on $\Sigma\mc{A}_X$. But then $\wt{\omega}$ will remain positive under small perturbations of $\omega$ and $J$. Similarly the condition of $\omega$ taming $J$ on $\ker \varphi$ is open. Consider $x \in X$ so that $\omega(v, J v) > 0$ for all $v \in \ker \varphi_x$. As $\wt{\omega}$ is continuous and $\Sigma\mc{A}_X$ is compact, there exists a neighbourhood $U$ of $\ker \varphi_x \cap \Sigma\mc{A}_X$ in $\Sigma\mc{A}_X$ on which $\wt{\omega}$ is positive. Points $x' \in X$ close to $x$ will then have $\ker \varphi_{x'} \subset U$ because $\ker \varphi$ is closed.
\ep
We will not use the associated notion of compatibility, where $J$ further leaves $\omega$ invariant, as we use \acs{}s as auxiliary structures to show non-degeneracy, and make use of the openness of this condition. For this reason, all \acs{}s will only be required to be continuous as this avoids arguments to ensure smoothness.

The following is the Lie algebroid version of \cite[Theorem 3.1]{Gompf04} and is our main tool to construct $\mc{A}_X$-\symp{}s; see \cite[Theorem 3.4]{CavalcantiKlaasse16} for the log-tangent version.
\begin{thm}\label{thm:lathurstontrick} Let $(\varphi,f)\colon (X,\mc{A}_X) \to (Y, \mc{A}_Y)$ be a Lie algebroid morphism between compact connected manifolds, $J$ an $\mc{A}_X$-\acs, $\omega_Y$ a $\mc{A}_Y$-symplectic form and $\eta$ a closed $\mc{A}_X$-two-form. Assume that
\bi
\item[i)] $\etwo{J}$ is $(\omega_Y,\varphi)$-tame;
\item[ii)] $\eta$ tames $J|_{\ker \varphi_x}$ for all $x \in X$.
\ei
Then $X$ admits an $\mc{A}_X$-\symp.
\end{thm}
The proof of this result is modelled on that by Gompf of \cite[Theorem 3.1]{Gompf04}.
\bp Let $t>0$ and consider the form $\omega_t := \varphi^* \omega_Y + t \eta$. We show $J$ is $\omega_t$-tame for $t$ small enough. By \autoref{prop:tamingopen} it is enough to show that there exists a $t_0 > 0$ so that $\omega_t(v, J v) > 0$ for every $t \in (0,t_0)$ and $v$ in the unit sphere bundle $\Sigma \mc{A}_X \subset \mc{A}_X$ with respect to some metric. For $v \in \mc{A}_X$ we have $\omega_t(v, J v) = \varphi^* \omega_Y(v,J v) + t \eta(v, J v)$. As $J$ is $(\omega_Y,\varphi)$-tame, the first term is positive for $v \in \mc{A}_X \setminus \ker \varphi$ and is zero otherwise. The second term $\eta(v, J v)$ is positive on $\ker \varphi$ because $J|_{\ker \varphi}$ is $\eta$-tame, hence is also positive for all $v$ in some neighbourhood $U$ of $\ker \varphi \cap \Sigma\mc{A}_X$ in $\mc{A}_X$ by openness of the taming condition. We conclude that $\omega_t(v, Jv) > 0$ for all $t > 0$ when $v \in U$. The function $v \mapsto \eta(v, J v)$ is bounded on the compact set $\Sigma\mc{A}_X \setminus U$. Furthermore, $\varphi^*\omega_Y(v, J v)$ is also bounded from below there by a positive constant, as it is positive away from $\ker \varphi$, and thus also away from $\ker \varphi \cap \Sigma\mc{A}_X \subset U$. But then $\omega_t(v, J v) > 0$ for all $0 < t < t_0$ for $t_0$ sufficiently small, so that $\omega_t$ is $\mc{A}_X$-symplectic for $t$ small enough.
\ep
Given a map $f\colon X \to Y$ we denote by $F_y = f^{-1}(y)$ for $y \in Y$ the level set, or \emph{fiber}, of $f$ over $y$. In order to meaningfully apply \autoref{thm:lathurstontrick} we must be able to construct closed $\mc{A}_X$-forms $\eta$ as in ii) of the statement. Note firstly that it suffices to have such closed forms in neighbourhoods of fibers, all lying in the same global Lie algebroid cohomology class.
\begin{prop}\label{prop:gluingfiberform} Assume that there exists a class $c \in H^2_{\mc{A}_X}(X;\R)$ and that for each $y \in Y$, $F_y$ has a neighbourhood $W_y$ with a closed $\mc{A}_X$-two-form $\eta_y \in \Omega^2_{\mc{A}_X}(W_y)$ such that $[\eta_y] = c|_{W_y} \in H^2_{\mc{A}_X}(W_y;\R)$, and $\eta_y$ tames $J|_{\ker \varphi_x}$ for all $x \in W_y$. Then there exists a closed $\mc{A}_X$-two-form $\eta$ such that $[\eta] = c$ and $\eta$ tames $J|_{\ker \varphi_x}$ for all $x \in X$.
\end{prop}
\bp Let $\xi \in \Omega^2_{\mc{A}_X}(X)$ be closed and such that $[\xi] = c$. Then for each $y \in Y$ we have $[\eta_y] = c|_{W_y} = [\xi]|_{W_y}$, so on $W_y$ we have $\eta_y = \xi + d_{\mc{A}_X} \alpha_y$ for some $\alpha_y \in \Omega^1_{{\mc{A}}_X}(W_y)$. As each $X \setminus W_y$ and hence $f(X \setminus W_y)$ is compact, each $y \in Y$ has a neighbourhood disjoint from $f(X\setminus W_y)$. Cover $Y$ by a finite amount of such open sets $U_i$ so that each $f^{-1}(U_i)$ is contained in some $W_{y_i}$. Let $\{\psi_i\}$ be a partition of unity subordinate to the cover $\{U_i\}$ of $Y$, so that $\{\psi_i \circ f\}$ is a partition of unity of $X$. Define an $\mc{A}_X$-two-form $\eta \in \Omega^2_{\mc{A}_X}(X)$ on $X$ via
\be
\eta := \xi + d_{\mc{A}_X}(\sum_i (\psi_i \circ f) \alpha_{y_i}) = \xi + \sum_i (\psi_i \circ f) d_{\mc{A}_X}\alpha_{y_i} + \sum_i \rho_{\mc{A}_X}^*(d\psi_i \circ df) \wedge \alpha_{y_i}.
\ee
Then $d_{\mc{A}_X}\eta = 0$ and $[\eta] = c$. As $\rho_{\mc{A}_X}^*(d\psi_i \circ df) = d\psi_i \circ df \circ \rho_{\mc{A}_X} = dg \circ \rho_{\mc{A}_Y} \circ \varphi$ using that $\varphi$ is a Lie algebroid morphism, the last of the above three terms vanishes when applied to a pair of vectors in $\ker \varphi_x$ for any $x \in X$, so on each $\ker \varphi_x$ we have
\be
\eta = \xi + \sum_i (\psi_i \circ f) d_{\mc{A}_X}\alpha_{y_i} = \sum_i (\psi_i \circ f)(\xi + d_{\mc{A}_X}\alpha_{y_i}) = \sum_i (\psi_i \circ f) \eta_{y_i}.
\ee
From the above we see that on $\ker \varphi$, the $\mc{A}_X$-form $\eta$ is a convex combination of $J$-taming $\mc{A}_X$-forms, so $J|_{\ker \varphi}$ is $\eta$-tame.
\ep
One can further look for (local or global) closed two-forms $\wt{\eta} \in \Omega^2(X)$ so that $\eta = \rho_{\mc{A}_X}^* \wt{\eta}$ satisfies hypothesis ii) of \autoref{thm:lathurstontrick}. When using such $\mc{A}_X$-forms which are pullbacks of regular forms, the behavior of the map $\rho_{\mc{A}_X}\colon \ker \varphi_x \to \ker df_x$ is important. We will see an example where there cannot be an $\eta$ of the form $\eta = \rho_{\mc{A}_X}^* \wt{\eta}$ making $\ker \varphi$ symplectic. This is unlike the situation for log-symplectic forms as studied in \cite{CavalcantiKlaasse16}, where the anchor of the log-tangent bundle provides an isomorphism between $\ker \varphi$ and $\ker df$.

Using \autoref{thm:lathurstontrick} we can prove the Lie algebroid version of Thurston's result for symplectic fiber bundles with two-dimensional fibers \cite{Thurston76}, adapting the proof by Gompf in \cite{Gompf01}.
\begin{thm}\label{thm:thurstonfibration} Let $(\varphi,f)\colon (X,\mc{A}_X^{2n}) \to (Y, \mc{A}_Y^{2n-2})$ be a Lie algebroid fibration between compact connected manifolds. Assume that $Y$ is $\mc{A}_Y$-symplectic and there exists a closed $\mc{A}_X$-two-form $\eta$ which is nondegenerate on $\ker \varphi$. Then $X$ admits an $\mc{A}_X$-\symp.
\end{thm}
\bp Let $\omega_Y \in {\rm Symp}(\mc{A}_Y)$ and choose $J_Y \in J(\omega_Y)$. Fix the orientation for $\ker \varphi$ so that $\eta$ is positive. Let $g$ be a metric on $\mc{A}_X$ and let $H \subset \mc{A}_X$ be the subbundle of orthogonal complements to $\ker \varphi$, so that $\varphi\colon H \to \mc{A}_Y$ is a fiberwise isomorphism. Define an $\mc{A}_X$-\acs $J$ by letting $J|_H = \varphi^* J_Y$, and on $\ker \varphi$, use the metric and define $J$ by $\frac{\pi}2$-counterclockwise rotation, demanding $\varphi$ is orientation preserving via the fiber-first convention. This determines $J$ uniquely on $\mc{A}_X$ by linearity. Moreover, $J$ is $(\omega_Y,\varphi)$-tame as $\varphi^* \omega_Y(v, J v) = \omega_Y(\varphi v, J_Y \varphi v) > 0$, for all $v \in \mc{A}_X \setminus \ker \varphi \cong H$. Further, $\eta$ tames $J$ on $\ker \varphi$ as $J$ is compatible with the orientation on $\ker \varphi$ determined by $\eta$. By \autoref{thm:lathurstontrick} we obtain an $\mc{A}_X$-\symp.
\ep
\begin{rem} One can combine \autoref{thm:thurstonfibration} with \autoref{prop:gluingfiberform} to obtain a statement requiring only the existence of local forms $\eta_y$ governed by a global cohomology class.
\end{rem}

Sadly for general Lie algebroid morphisms there is no direct analogue of homological essentialness of generic fibers to replace the hypothesis of the existence of $\eta$ as in \autoref{thm:lathurstontrick} and \autoref{thm:thurstonfibration}. In other words, we cannot replace the condition on the existence of $\eta$ by demanding that the generic fiber $F$ satisfies $[F] \neq 0 \in H_2(X;\R)$. This is a well-known necessary condition when constructing symplectic-like structures out of Lefschetz-type fibrations using Gompf--Thurston techniques. The reason for the lack of such an analogue is again the behavior of the map $\rho_{\mc{A}_X}\colon \ker \varphi \to \ker df$. While the codomain can be seen as the tangent space to the fiber at regular points, the domain cannot. For Lie algebroid submersions, surjectivity is only demanded of $\varphi$, not of $df$, and there is no duality pairing between homology and Lie algebroid cohomology in general.
\begin{rem} In the special cases of the log-tangent bundle and the elliptic tangent bundle one can state such an analogue (see Theorems \ref{thm:lfibrationsgcs} and \ref{thm:fibrationsgcs}, and \cite{CavalcantiKlaasse16}), precisely because we understand their respective Lie algebroid cohomologies (see \autoref{thm:loglacohomology} and \autoref{thm:elllacohomology}).
\end{rem}

Recall that we defined a Lie algebroid \lf to be a Lie algebroid fibration away from the set of Lefschetz singularities, which must be a subset of the isomorphism locus of $\mc{A}_X$. Note that $X_{\mc{A}_X}$ is either empty or of codimension zero, so that we can restrict $\mc{A}_X$ to $X_{\mc{A}_X}$ and again obtain a Lie algebroid. In this way the inclusion $i\colon X_{\mc{A}_X} \hookrightarrow X$ covers a Lie algebroid morphism $\iota$, and we obtain a restriction map in cohomology $\iota^*\colon H_{\mc{A}_X}(X) \to H_{\mc{A}_X}(X_{\mc{A}_X})$. The anchor gives a bijection $\rho_{\mc{A}_X}^*\colon \Omega^\bullet(X_{\mc{A}_X}) \to \Omega^\bullet_{\mc{A}_X}(X_{\mc{A}_X})$ and hence an isomorphism $\rho_{\mc{A}_X}^*\colon H_{\rm dR}(X_{\mc{A}_X}) \to H_{\mc{A}_X}(X_{\mc{A}_X})$ in cohomology over $X_{\mc{A}_X}$. 

Assuming ${\rm rank}(\mc{A}_X) = {\rm rank}(\mc{A}_Y) + 2$, the second hypothesis in \autoref{thm:lathurstontrick} demands the existence of a class $c \in H_{\mc{A}_X}(X)$ such that the de Rham class $(\rho_{\mc{A}_X}^*)^{-1} \iota^*(c)$ evaluates nonzero on each generic fiber. Even if such a class exists, a similar statement must hold over $X \setminus X_{\mc{A}_X}$. Hence, we need the existence of a two-form $\eta \in \Omega^2_{\rm cl}(X;\mc{A}_X)$ such that $(\rho_{\mc{A}_X}^*)^{-1} \iota^* [\eta]$ evaluates nonzero on the fibers, and $\eta|_{\ker \varphi}$ is nowhere zero over $X \setminus X_{\mc{A}_X}$.

The notion of a Lie algebroid \lf is such that the hypotheses of \autoref{thm:lathurstontrick} are still satisfied, when in dimension four.
\begin{thm}\label{thm:thurstonlfibration} Let $(\varphi,f)\colon \mc{A}^4_{X} \to \mc{A}^2_{\Sigma}$ be a Lie algebroid \lf with connected fibers between compact connected manifolds. Assume that $\mc{A}_\Sigma$ admits a symplectic structure and there exists a closed $\mc{A}_X$-two-form $\eta$ such that $(\rho_{\mc{A}_X}^*)^{-1} \iota^* [\eta]$ evaluates nonzero on the fibers, and $\eta|_{\ker \varphi}$ is nondegenerate over $X \setminus X_{\mc{A}_X}$. Then $X$ admits an $\mc{A}_X$-\symp.
\end{thm}
The proof of this result is similar to \cite[Theorem 4.7]{CavalcantiKlaasse16}, which in turn is based on Gompf's proof of \cite[Theorem 10.2.18]{GompfStipsicz99} using almost-complex structures. We will use \autoref{prop:gluingfiberform} to glue together the given form $\eta$ with a suitable adaptation in the isomorphism locus of $\mc{A}_X$.
\begin{rem} Note that if $\mc{A}_X = TX$, then $X_{\mc{A}_X} = X$ and we are demanding that the generic fiber is homologically essential, recovering the result by Gompf \cite[Theorem 10.2.18]{GompfStipsicz99}.
\end{rem}
\bp If $X_{\mc{A}_X}$ is empty, $\eta|_{\ker \varphi}$ is nondegenerate everywhere and $\varphi$ is a Lie algebroid fibration, so that the result follows from \autoref{thm:thurstonfibration}. If not, then $X$ is four-dimensional. Denote $\xi := (\rho_{\mc{A}_X}^*)^{-1} \iota^* \eta \in \Omega^2(X_{\mc{A}_X})$ and $c = [\xi] \in H^2_{{\rm dR}}(X_{\mc{A}_X})$. Note that $\xi$ orients the generic fibers $F$, as these are two dimensional. Recall that $\Delta \subset X_{\mc{A}_X}$ is the set of Lefschetz singularities, and let $\Delta' := f(\Delta)$ be the set of singular values of $f|_{X_{\mc{A}_X}}$.

Let $\omega_\Sigma$ be an $\mc{A}_\Sigma$-\symp and use the proof of \autoref{thm:thurstonfibration} to obtain a $(\varphi,\omega_\Sigma)$-tame almost-complex structure $J$ on $\mc{A}_X$ over $X \setminus f^{-1}(\Delta')$ compatible with the orientation on $\ker \varphi$, noting that $f\colon X \setminus f^{-1}(\Delta') \to \Sigma \setminus \Delta'$ is a Lie algebroid fibration. As $\Delta'$ is contained in $\Sigma_{\mc{A}_\Sigma}$, let $V \subseteq \Sigma$ be the disjoint union of open balls $V_y$ disjoint from $Z_\Sigma$ and centered at each point $y \in \Delta'$. Let $W := f^{-1}(V) \subset X$ be the union of the neighbourhoods $W_y := f^{-1}(V_y)$ of singular fibers $F_y$. Let $C \subset X$ be the disjoint union of open balls $C_y \subseteq \Sigma_{\mc{A}_\Sigma}$ centered at each point $f^{-1}(y) = x$ for all $y \in \Delta'$ so that on each ball $f$ takes on the local form in \autoref{defn:lefschetzfibration}. Possibly shrink $C$ so that $\overline{C_y} \subset W_y$. The local description of $f$ gives an \acs on $C$ with the fibers being holomorphic, and we glue this to the existing almost-complex structure $J$ on $X \setminus C$. This gives a global $(\omega_\Sigma, \varphi)$-tame $\mc{A}_X$-\acs $J$.

As $\eta|_{\ker \varphi}$ is nondegenerate over $X \setminus X_{\mc{A}_X}$, the same is true in a neighbourhood $S$ around $X \setminus X_{\mc{A}_X}$ disjoint from $\Delta$. Let $y \in Y \setminus f(S)$ be given. If $y \not\in \Delta'$, let $D_y \subset Y \setminus f(S)$ disjoint from $\Delta'$ be a disk containing $y$, fully contained in a trivializing neighbourhood of $f$ around $y$. Define $W_y := f^{-1}(D_y) \cong D_y \times F_y$, with projection map $p:\, W_y \to D_y$. Let $\eta_y'$ be an area form on $F_y$ compatible with the preimage orientation. Define $\eta_y = \lambda_y \rho_X^* p^* \eta_y'$, where $\lambda_y \in \R$ is chosen such that $\langle [F], c \rangle = \langle [F], \eta_y \rangle$. As $H_2(W_y;\R)$ is generated by $[F_y]$, it follows that $[\eta_y] = c|_{W_y} \in H^2(W_y;\R)$. But then $\eta_y$ tames $J$ on $\ker \varphi \cong \ker df$ for all $x \in W_y$, as the restriction of $(\rho_X^*)^{-1} \eta_y = \lambda_y p^* \eta_y'$ is an area form for that fiber.

If $y \in \Delta'$, the singular fiber $F_y$ either is indecomposable or consists of two irreducible components $F_y^\pm$ which satisfy $[F_y^+] \cdot [F_y^-] = 1$ and $[F_y^\pm]^2 = -1$, see \cite{GompfStipsicz99}. In the latter case, note that $0 < 1 = \langle c, [F] \rangle = \langle c, [F_y]\rangle = \langle c, [F_y^+]\rangle + \langle c, [F_y^-]\rangle$. If either term is nonpositive assume without loss of generality that $\langle c, [F_y^-]\rangle = r \leq 0$. Define $c' := c + (\frac12 - r) c_y^+$, where $c_y^+ \in H^2(X;\R)$ is a class dual to $[F_y^+]$. As $[F_y] \cdot [F_y^\pm] = 0$ we then have $\langle c', [F] \rangle = \langle c, [F] \rangle > 0$, and furthermore $\langle c', [F_y^+]\rangle = \langle c, [F_y^+] \rangle +(\frac12 - r) > 0$ and $\langle c', [F_y^-]\rangle = \frac12 > 0$. Moreover, as different fibers do not intersect, we have $c|_{W_{y'}} = c'|_{W_{y'}}$ for $y' \neq y$. After finitely many repetitions, at most once for each $y \in \Delta'$, one obtains a class, again denoted by $c$, pairing positively with every fiber component (see \cite[Exercise 10.2.19]{GompfStipsicz99}).

Return to $y \in \Delta'$ and let $\sigma$ be the standard symplectic form on $C_y$ given locally in real coordinates by $\sigma = dx_1 \wedge dy_1 + dx_2 \wedge dy_2$, where $z_i = x_i + i y_i$. As all fibers $F_y'$ in $C_y$ are holomorphic, $\rho_X^* \sigma|_{F_{y'} \cap C_y}$ tames $J$ for all $y' \in f(C_y)$, so that $\rho_X^* \sigma$ tames $\btwo{J}$ on $C_y$. Let $\sigma_y$ be an extension of $\sigma$ to $F_y$ as a positive area form with total area $\langle \sigma_y, [F_y] \rangle$ equal to $\langle c, [F_y] \rangle$. Let $p\colon W_y \to F_y$ be a retraction and let $f\colon C_y \to [0,1]$ be a smooth radial function so that $f \equiv 0$ in a neighbourhood of $x = f^{-1}(y) \in \Delta$ and $f \equiv 1$ in a neighbourhood of $\partial C_y$, which is smoothly extended to $W_y$ by being identically $1$ outside $C_y$. On the ball $C_y$, the form $\sigma$ is exact, say equal to $\sigma = d\alpha$ for $\alpha \in \Omega^1(C_y)$. Define a two-form $\eta_y'$ on $W_y$ by $\eta_y' := p^*(f\sigma_y) + d((1-f)\alpha)$, which is closed as $f\sigma_y$ is a closed area form on $F_y$. Near $x$ we have $f \equiv 0$ so that $\eta_y' = d\alpha = \sigma$. Set $\eta_y := \rho_X^* \eta_y'$. Then there $\eta_y = \rho_X^* \sigma$ tames $J$, hence in particular tames $J|_{\ker \varphi}$. Similarly, $\sigma_y$ is an area form on $F_y \setminus \{x\}$ for the orientation given by $J$ under the isomorphism by $\rho_X$. But then $\rho_X^* \sigma_y$ tames $J$ on $\ker \varphi_y \cong \ker df = TF_y$ on $F_y$, so that the same holds for $\eta_y$ as this condition is convex. By openness of the taming condition, shrinking $V_y$ and hence $W_y$ and possibly $C_y$ we can ensure that $\eta_y$ tames $J|_{\ker \varphi}$ on $W_y$. Finally, note that $[\eta_y'] = c|_{W_y} \in H^2(W_y;\R)$ by construction. For any point $y \in f(S)$, take the Lie algebroid two-form $\eta_y := \eta$ on the neigbourhood $W_y := S$ of $F_y$. We have now obtained the required neighbourhoods $W_y$ and forms $\eta_y$ for all $y \in Y$ to apply \autoref{prop:gluingfiberform} and obtain a Lie algebroid closed two-form again denoted by $\eta$ such that $\eta$ tames $J$ on $\ker \varphi_x$ for all $x \in X$. By \autoref{thm:lathurstontrick} we obtain an $\mc{A}_X$-symplectic structure on $X$.
\ep
\begin{rem} When both $\mc{A}_X$ and $\mc{A}_\Sigma$ are log-tangent bundles, \autoref{thm:thurstonlfibration} recovers \cite[Theorem 4.7]{CavalcantiKlaasse16}, upon combination with the fact that $\ker \varphi$ and $\ker df$ are pointwise isomorphic.
\end{rem}

\section{Boundary maps and \placeholder{}s}
\label{sec:boundarymaps}
In this section we introduce the notion of a boundary map, which is a map degenerating suitably on a submanifold. When this submanifold has codimension two these will supply us with morphisms from elliptic to log divisors, and hence with Lie algebroid morphisms from the respective elliptic to log-tangent bundles by \autoref{prop:elltologmorphism}. After this we define the notion of a \placeholder, which is a Lefschetz-type fibration that can be interpreted as a Lie algebroid \lf between these Lie algebroids.
\subsection{The normal Hessian}
Let $(X,D)$ be a \emph{pair}, i.e.\ $X$ is a manifold and $D \subseteq X$ is a submanifold.
\begin{defn} A \emph{map of pairs} $f\colon (X,D) \to (Y,Z)$ is a map $f\colon X \to Y$ such that $f(D) \subseteq Z$. A \emph{strong map of pairs} is a map of pairs $f\colon (X,D) \to (Y,Z)$ such that $f^{-1}(Z) = D$.
\end{defn}
Given a map of pairs $f\colon (X,D) \to (Y,Z)$, note that $df(TD) \subseteq TZ$, so that $df$ induces a map $\nu(df)\colon ND \to NZ$ between normal bundles. When $f$ is a strong map of pairs and $f$ is transverse to $Z$, the map $\nu(df$) is an isomorphism.

Consider a map of pairs $f\colon (X,D) \to (Y,Z)$ such that $\nu(df)\colon ND \to NZ$ is the zero map. Equivalently, one can assume that ${\rm im}\, df \subset TZ$. Let $z_1, \dots, z_\ell$ be local defining functions for $Z$ and consider their pullbacks $h_i := f^*(z_i)$ for $i = 1, \dots, \ell$. As $f(D) \subseteq Z$, the functions $h_i$ vanish on $D$. Moreover, because $\nu(df)$ is the zero map, the derivatives $dh_i$ vanish on $D$ as well. Consequently, we can consider their Hessians $H(h_i)\colon {\rm Sym}^2(TX) \to \R$, which descend to maps $H(h_i)\colon {\rm Sym}^2(ND) \to \R$. As the differentials $dz_i$ span $N^*Z$, these combine to give a map $H^\nu(f)\colon {\rm Sym}^2(ND) \to f^*(NZ)$, which one checks to be invariantly defined.
\begin{defn} Let $f\colon (X,D) \to (Y,Z)$ such that ${\rm im}\, df \subset TZ$. The \emph{normal Hessian} of $f$ along $D$ is the map $H^\nu(f): {\rm Sym}^2(ND) \to f^*(NZ)$ over $D$.
\end{defn}
When ${\rm codim}(Z) = 1$, the normal Hessian $H^\nu(f)$ can be viewed as the matrix of second partial derivatives of the coordinate function of $f$ normal to $Z$ in directions normal to $D$.
\subsection{Boundary maps}
Let $f\colon (X,D) \to (Y,Z)$ be a strong map of pairs with $Z$ a hypersurface and ${\rm codim}(D) \geq 2$. Then $f$ cannot be transverse to $Z$, as then $f^{-1}(Z) = D$ would be of codimension one. As $Z$ is of codimension one, $f$ being transverse to $Z$ is equivalent to $\nu(df)\colon ND \to NZ$ being nonzero, hence in this case $\nu(df)\colon ND \to NZ$ is the zero map. Equivalently we have ${\rm im}\, df \subset TZ$, so that the normal Hessian of $f$ along $D$ is well-defined.
\begin{defn} Let $f\colon (X,D) \to (Y,Z)$ be a strong map of pairs, $Z$ a hypersurface, and ${\rm codim}(D) \geq 2$. Then $f$ is a \emph{boundary map} if its normal Hessian $H^\nu(f)$ is definite along $D$.
\end{defn}
As $Z$ is a hypersurface, $NZ$ is one-dimensional. Because of this, $H^\nu(f)$ being definite makes sense, as locally it is a map $H^\nu(f)\colon {\rm Sym}^2(\R^d) \to \R$ where $d = {\rm codim}(D)$. The choice of the name will become clearer after establishing some properties (see \autoref{rem:bdrymapname}). We will also call $f$ a \emph{codimension-$k$ boundary map} if ${\rm codim}(D) = k$, and sometimes implicitly assume that ${\rm codim}(D) = 2$ when $\dim(X) = 4$. Indeed, the main reason for introducing the notion of a boundary map comes from \autoref{prop:bdrymapmorphdivisors} below, where ${\rm codim}(D) = 2$.
\begin{rem} Let $f:\, (X,D) \to (Y,Z)$ be a codimension-$k$ boundary map and $f'':\, (X',D') \to (X,D)$ and $f':\, (Y,Z) \to (Y',Z')$ strong maps of pairs with $f''$ transverse to $D$ and $f'$ transverse to $Z'$. Then $f'' \circ f \circ f':\, (X',D') \to (Y',Z')$ is a codimension-$k$ boundary map.
\end{rem}
\begin{rem} Assuming that ${\rm codim}(D) \geq 2$ is only done to ensure that $\nu(df)$ is the zero map, as is required for the definition of the normal Hessian. If ${\rm codim}(D) = 1$ yet this condition holds, it makes sense to talk about codimension-one boundary maps.
\end{rem}
Codimension-two boundary maps naturally give rise to morphisms of divisors between elliptic and log divisors respectively, and hence to Lie algebroid morphisms by \autoref{prop:elltologmorphism}.
\begin{prop}\label{prop:bdrymapmorphdivisors} Let $f\colon (X,D) \to (Y,Z)$ be a map of pairs with $Z$ a hypersurface and ${\rm codim}(D) = 2$. Then $f$ is a boundary map if and only if $I_D := f^* I_Z$ is an elliptic ideal and $f$ a morphism of divisors.
\end{prop}
In other words, a codimension-two boundary map uniquely specifies a compatible elliptic divisor structure on $D$. 
\bp Assume that $f$ is a boundary map and consider $I_D = f^* I_Z$. Let $z$ be a local defining function for $Z$ so that locally $I_Z = \langle z \rangle$, and hence $I_D = \langle f^*(z) \rangle$. As $f$ is a boundary map, $H^\nu(f)$ is definite, so that $f^*z$ specifies the germ of a definite Morse--Bott function around $D$. By the discussion above \autoref{prop:ellmorsebott} we see that $I_D$ is an elliptic ideal specifying an elliptic divisor structure on $D$ by \autoref{prop:locprincideal}. The map $f$ is a morphism of divisors by construction. Alternatively, one shows that for $(L,s)$ the log divisor determined by $Z$, the pair $(R,q) := (f^* L, f^* s)$ is an elliptic divisor with $D_q = D$. The converse is similar, using again that $f$ has definite normal Hessian if and only if $f^*(z)$ is locally Morse--Bott of even index around $D$, where $z$ is a local defining function for $Z$.
\ep
\begin{cor}\label{cor:bdrymaplamorphism} Let $f\colon (X,D) \to (Y,Z)$ be a codimension-two boundary map. Then $df$ induces a Lie algebroid morphism $(\varphi,f)\colon TX(-\log |D|) \to TY(- \log Z)$ for the divisor structures of \autoref{prop:bdrymapmorphdivisors}.
\end{cor}
The pointwise conclusion of the Morse--Bott lemma, \autoref{lem:morsebott}, provides a local form for boundary maps around points in $D$.
\begin{prop}\label{prop:bdrymaplocalcoord} Let $f\colon (X^n,D^{k}) \to (Y^m,Z^{m-1})$ be a boundary map and $x \in D$. Then around $x$ and $f(x) \in Z$ there exist coordinates $(x_1,\dots,x_n)$ and $(z,y_2,\dots,y_m)$ for which ${\{x_1 = \dots = x_{n-k} = 0\} = D}$ and $\{z = 0\} = Z$ such that for some map $g\colon \R^n \to \R^{k-1}$ we have $f(x_1,\dots,x_n) = (x_1^2 + \dots + x_{n-k}^2, g(x_1,\dots,x_n))$.
\end{prop}
\bp Let $x \in D$ and $z$ be a local defining function for $Z$ around $f(x) \in Z$. As $f$ is a boundary map, the proof of \autoref{prop:bdrymapmorphdivisors} shows that $f^*(z)$ is a local Morse--Bott function of index zero around $x$, after possibly replacing $z$ by $-z$. By \autoref{lem:morsebott}, after trivializing $ND$, there exist coordinates $(x_1,\dots,x_n)$ of $X$ around $x$ with $\{x_1 = \dots = x_{n-k} = 0\} = D$ such that $f^*(z)(x_1,\dots,x_n) = x_1^2 + \dots + x_{n-k}^2$. Complete $z$ to a coordinate system $(z,y_2,\dots,y_m)$ of $Y$ around $f(x)$. Then $f(x_1,\dots,x_n) = (f^*(z)(x_1,\dots,x_n),g(x_1,\dots,x_n)) = (x_1^2 + \dots + x_{n-k}^2, g(x_1,\dots,x_n))$ for some $g\colon \R^n \to \R^{k-1}$.
\ep
Using either \autoref{prop:bdrymaplocalcoord} or the proof \autoref{prop:bdrymapmorphdivisors}, the Lie algebroid morphisms $(\varphi,f)\colon TX(-\log|D|) \to TY(-\log Z)$ obtained from boundary maps $f\colon (X,D) \to (Y,Z)$ using \autoref{cor:bdrymaplamorphism} have the following extra property: for any $k$-form $\alpha \in \Omega^k(\log Z)$, we have ${\rm Res}_\theta(\varphi^* \alpha) = 0$, because ${\rm Res}_\theta f^*(z) = 0$ for any local defining function $z$ for $Z$.
\subsection{Fibrating boundary maps}
We next introduce specific boundary maps by demanding submersiveness on $D$.
\begin{defn} A \emph{fibrating boundary map} is a boundary map $f\colon (X,D) \to (Y,Z)$ such that $f|_{D}\colon D \to Z$ is submersive.
\end{defn}
Note that it is not required that $D$ surjects onto $Z$. For fibrating boundary maps we can improve upon \autoref{prop:bdrymaplocalcoord}, obtaining again a local form around points in $D$.
\begin{prop}\label{prop:fibratingpointwise} Let $f\colon (X^n,D^{k}) \to (Y^{m},Z^{m-1})$ be a fibrating boundary map. Then around points $x \in D$ and $f(x) \in Z$ there exist coordinates $(x_1,\dots,x_n)$ and $(z,y_2,\dots,y_m)$ with $\{x_1 = \dots = x_{n-k} = 0\} = D$ and $\{z = 0\} = Z$ such that $f(x_1,\dots,x_n) = (x_1^2 + \dots + x_{n-k}^2, x_{n-m+1}, \dots, x_{n})$.
\end{prop}
In other words, we can simultaneously put both the components of $f$ normal and tangent to $D$ in standard form, and obtain  a commuting diagram of normal bundles near $x$ and $f(x)$:
\begin{center}
	\begin{tikzpicture}
	\matrix (m) [matrix of math nodes, row sep=2.5em, column sep=2.5em,text height=1.5ex, text depth=0.25ex]
	{	 ND & NZ \\ D & Z\\};
	\path[-stealth]
	(m-1-1) edge node [above] {$f$} (m-1-2)
	(m-1-1) edge node [left] {${\rm pr}_D$} (m-2-1)
	(m-1-2) edge node [right] {${\rm pr}_Z$} (m-2-2)
	(m-2-1) edge node [above] {$f|_D$} (m-2-2);
	\end{tikzpicture}
\end{center}
Call a finite collection of functions an \emph{independent set} at $p$ if their differentials are everywhere linearly at $p$. By the implicit function theorem, an independent set can be completed to a coordinate system in a neighbouhood of $p$. Independence is preserved under pulling back along a submersion.
\bp Choose a tubular neighbourhood embedding $\Phi\colon NZ \to V \subset Y$ for ${\rm pr}_Z\colon NZ \to Z$ and let $z$ be a local defining function for $Z$ on an open subset $V' \subset V$ of $f(x)$. Let $U := f^{-1}(V') \subset X$. Choose coordinate functions $y_2,\dots,y_m\colon V' \to \R$ for $Z$. Then $\{y_2,\dots,y_m\}$ forms an independent set on $Z$, and because ${\rm pr}_Z$ is a submersion, the same is true for $\{z, {\rm pr}_Z^*(y_2), \dots, {\rm pr}_Z^*(y_m)\}$ on $Y$. By \autoref{prop:bdrymaplocalcoord}, after possibly shrinking $U$ and $V'$ and changing $z$ to $-z$, there exist tubular neighbourhood coordinates $(x_1,\dots,x_n)$ around $x$ such $f^*(z) = x_1^2 + \dots + x_{n-k}^2$. Consider $\{x_1,\dots,x_{n-k},f^* {\rm pr}_Z^*(y_j)\}$, which is an independent set on $X$, using submersiveness of $f|_D$. Complete this to a coordinate system $\{x_1,\dots,x_{n-m+1},f^* {\rm pr}_Z^*(y_j)\}$ on $X$, and relabel these as $(x_1,\dots,x_n)$. Using these coordinates on $X$ and the coordinates $(z,\pi_Z^*(y_j))$ on $Y$, the map $f$ is given by $f(x_1,\dots,x_n) = (x_1^2 + \dots + x_{n-k}^2, x_{n-m+1}, \dots, x_{n})$ as desired.
\ep
The normal form result of \autoref{prop:fibratingpointwise} immediately implies the following.
\begin{cor}\label{cor:fibratingsubm} Let $f\colon (X,D) \to (Y,Z)$ be a fibrating boundary map. Then $f$ is submersive in a punctured neighbourhood around $D$.
\end{cor}
Consequently, fibrating boundary maps have well-defined fibers of dimension $\dim(X) - \dim(Y)$ near $D$, and of dimension $\dim(X) - \dim(Y) - {\rm codim}(D) + 1$ on $D$. In particular, when $f$ is a fibrating boundary map, $D$ will be a fiber bundle over certain components of $Z$.
\begin{rem}\label{rem:fibrliealgdmorph} An alternative way of viewing the proof of \autoref{cor:fibratingsubm} when ${\rm codim}(D) = 2$ uses \autoref{prop:bdrymapmorphdivisors}. Namely, let $f:\, (X,D) \to (Y,Z)$ give rise to a Lie algebroid morphism $(\varphi,f)\colon TX(-\log |D|) \to TY(-\log Z)$. While $df\colon ND \to NZ$ is the zero map, $\varphi|_D$ is surjective because $f$ is fibrating. This is an open condition, so that $\varphi$ is surjective in a neighbourhood around $D$. On $X \setminus D$, the Lie algebroid $TX(-\log |D|)$ is isomorphic to $TX$, hence $f$ is submersive there. We see that fibrating boundary maps lead to Lie algebroid morphisms which are fibrations near $D$.
\end{rem}
We proceed to obtain similar normal bundle commutativity around components of $D$ whose image is coorientable. To prove this, we use a result by Bursztyn--Lima--Meinrenken \cite{BursztynLimaMeinrenken16} on normal bundle embeddings, which we now describe. Let $M \subseteq X$ be a submanifold and $NM$ its normal bundle. Denote by $\mc{E}_M \in \mf{X}(NM)$ the associated Euler vector field and given $V \in \mf{X}(X)$ tangent to $M$, let $\nu(X) \in \mf{X}(NM)$ be its linear approximation.
\begin{defn} Let $M \subseteq X$ be a submanifold. A vector field $V \in \mf{X}(X)$ is \emph{Euler-like} along $M$ if it is complete, and satisfies $V|_M = 0$ and $\nu(V) = \mc{E}_M$.
\end{defn}
A \emph{strong tubular neighbourhood embedding} for $M \subseteq X$ is an embedding $\Phi\colon NM \to X$ taking the zero section of $NM$ to $M$, and such that the linear approximation $\nu(\Phi)\colon (NM,M) \to (X,M)$ is the identity map. The following proposition says that Euler-like vector fields give rise to strong tubular neighbourhood embeddings.
\begin{prop}[{\cite[Proposition 2.6]{BursztynLimaMeinrenken16}}]\label{prop:strongtubnhood} Let $M \subseteq X$ be a submanifold and $V \in \mf{X}(X)$ Euler-like along $M$. Then there exists a unique strong tubular neighbourhood embedding $\Phi\colon NM \to X$ such that $\Phi_*(\mc{E}_M) = V$.
\end{prop}
We use this result to construct compatible tubular neighbourhood embeddings for fibrating boundary maps around components of $D$ whose image is coorientable.
\begin{prop}\label{prop:fibratingsemiglobal} Let $f\colon (X,D) \to (Y,Z)$ be a fibrating boundary map and $D_j \subseteq D$ a connected component such that $f(D_j) =: Z_j \subset Z$ is coorientable. Then there exist a defining function $z$ for $Z_j$ and tubular neighbourhood embeddings $\Phi_{D_j}\colon \wt{U} \to U \subset X$ for $D_j$ and $\Phi_{Z_j}\colon\wt{V} \to V \subset Y$ for $Z_j$ such that $\Phi_{D_j}^* f^*(z) = Q_{f^*(z)} \in \Gamma(D_j; {\rm Sym}^2 N^* D_j)$ and the following diagram commutes:
\begin{center}
	\begin{tikzpicture}
	\matrix (m) [matrix of math nodes, row sep=2.5em, column sep=2.5em,text height=1.5ex, text depth=0.25ex]
	{	 \wt{U} & U & V & \wt{V}\\ D_j & D_j & Z_j & Z_j\\};
	\path[-stealth]
	(m-1-1) edge node [above] {$\Phi_{D_j}$} (m-1-2)
	(m-1-1) edge node [left] {${\rm pr}_{D_j}$} (m-2-1)
	(m-2-1) edge node [above] {${\rm id}$} (m-2-2)
	(m-1-2) edge node [above] {$f$} (m-1-3)
	(m-2-2) edge node [above] {$f|_D$} (m-2-3)
	(m-1-4) edge node [above] {$\Phi_{Z_j}$} (m-1-3)
	(m-2-4) edge node [above] {${\rm id}$} (m-2-3)
	(m-1-4) edge node [right] {${\rm pr}_{Z_j}$} (m-2-4);
	
	\draw[right hook->]
	(m-2-2) edge (m-1-2)
	(m-2-3) edge (m-1-3);
	\end{tikzpicture}
\end{center}
\end{prop}
As a consequence, when $Z$ is separating this result implies a normal form for $f$ around any point in $Z$ and its entire inverse image, as then a global defining function for $Z$ can be used.
\bp Let $z\colon V \to \R$ be a defining function for $Z_j$ and let $U' \subset f^{-1}(V)$ be the connected component containing $D_j$. Using \autoref{lem:morsebott} applied to $f^*(z)$, shrink $U'$ so that $U' = \Phi_{D_j}(\wt{U}')$ for some tubular neighbourhood embedding $\Phi_{D_j}\colon \wt{U}' \to U$ of $D_j$. Choose a tubular neighbourhood embedding $\Phi_{Z_j}\colon \wt{V} \to V$ for $Z_j$. For $x \in D_j$, use the proof of \autoref{prop:fibratingpointwise} (possibly changing $z$ to $-z$) to obtain an open $U_x \subset U'$ containing $x$ and coordinates $(x_1,\dots,x_n)$ so that $f^*(z) = x_1^2 + \dots + x_{n-k}^2$. Note that $U'$ is connected so that $f^*(z)$ has a fixed sign. Consider $U := \cup_{x \in D_j} U_x \subset U'$ and extract a finite subcover $\{U_\alpha\}_{\alpha \in I}$. On each set $U_\alpha$, define the vector field $\mc{E}_\alpha := x_1 \partial_{x_1} + \dots + x_{n-k} \partial_{x_{n-k}}$. It satisfies $\mc{L}_{\mc{E}_\alpha} f^*(z) = 2 f^*(z)$, and $f_* \Phi_{\alpha *} \mc{E}_\alpha = z \partial_z = \mc{E}_{Z_j}$. Let $\{\psi_\alpha\}_{\alpha \in I}$ be a partition of unity subordinate to $\{U_\alpha\}_{\alpha \in I}$ and define $\mc{E} := \sum_{\alpha \in I} \psi_{\alpha} \mc{E}_\alpha$ on $U'$. Then $\mc{E}|_{D_j} = \sum_{\alpha \in I} \psi_\alpha \mc{E}_\alpha|_{D_j} = \sum_{\alpha \in I} \psi_\alpha \cdot 0 = 0$, and $\nu(\mc{E}) = \sum_{\alpha \in I} \psi_\alpha \nu(\mc{E}_\alpha) = \sum_{\alpha \in I} \psi_{\alpha} \mc{E}_{D_j} = \mc{E}_{D_j}$. Ensure that $\mc{E}$ is complete by multiplying by a bump function and shrinking $U'$, so that $\mc{E} \in \mf{X}(U)$ is Euler-like along $D_j$. Moreover, $\mc{L}_{\mc{E}} f^*(z) = \sum_{\alpha \in I} \mc{L}_{\psi_\alpha \mc{E}_\alpha} f^*(z) = \sum_{\alpha \in I} \psi_\alpha \mc{L}_{\mc{E}_\alpha} f^*(z) = \sum_{\alpha \in I} \psi_\alpha 2 f^*(z) = 2 f^*(z)$. Finally, $f_* \mc{E} = \sum_{\alpha \in I} f_* (\psi_\alpha \mc{E}_\alpha) = \sum_{\alpha \in I} \psi_\alpha \mc{E}_Z = \mc{E}_Z$. Use \autoref{prop:strongtubnhood} and possibly shrink $U$ to obtain a tubular neighbourhood embedding $\Phi_{D_j}'\colon \wt{U} \to U$ such that $\Phi_* \mc{E}_{D_j} = \mc{E}$, $f_* \mc{E} = \mc{E}_Z$, and $\mc{L}_{\mc{E}} f^*(z) = 2 f^*(z)$. This embedding satisfies all desired properties. Moreover, $\Phi^{'*}_{D_j} f^*(z) = Q_{f^*(z)} \in \Gamma(D_j; {\rm Sym}^2(N^*D_j))$ as $g :=\Phi^{'*}_{D_j} f^*(z)$ is smooth and satisfies $\mc{L}_{\mc{E}} g = 2 g$.
\ep
We can now obtain topological information about the generic fibers of $f$ near $D$. The following result is the main reason for wanting compatible tubular neighbourhood embeddings.
\begin{prop}\label{prop:fiberbundle} Let $f\colon (X^n,D^{k}) \to (Y^m,Z^{m-1})$ be a fibrating boundary map with $Z$ separating. Denote the fiber of $f|_D\colon D \to Z$ by $F_D^{k-m+1}$. Then the fiber $F^{n-m}$ of $f$ near $D$ is an $S^{n-k-1}$-sphere bundle over $F_D$.
\end{prop}
\bp Apply \autoref{prop:fibratingsemiglobal} to each connected component of $D$ using a single global defining function $z$ for $Z$. Let $y \in \wt{V} \setminus Z$. As $f$ is submersive on $\wt{U} \setminus D$, consider $F_y = \Phi_{D}^{-1} \circ f^{-1} \circ \Phi_{Z}(y)$, which is of dimension $n-m$. Consider ${\rm pr}_Z(y) \in Z$ and its $(k-m+1)$-dimensional fiber $F_{D,{\rm pr}_Z(y)} = f|_D^{-1}({\rm pr}_Z(y))$. Because the tubular neighbourhood embeddings are compatible with $f$, we have $F_y = {\rm pr}_D^{-1}(F_{D,{\rm pr}_Z(y)})$, noting that ${\rm pr}_D$ is submersive. As a point in $\wt{V} \subset NZ$ is given by a point in $Z$ together with a distance, the fiber of ${\rm pr}_D\colon F_y \to F_{D,{\rm pr}_Z(y)}$ is given by a sphere, consisting of all points with fixed distance above the corresponding point in $D$.
\ep
\begin{cor}\label{cor:fibertori} Let $f\colon (X^4,D^2) \to (Y^2,Z^1)$ be a codimension-two fibrating boundary map. Then the generic fibers of $f$ near $D$ are unions of tori.
\end{cor}
\bp This follows from \autoref{prop:fiberbundle} using $n = 4$ and $k = m = 2$ by a dimension count. The generic fibers $F$ of $f$ near $D$ satisfy ${\rm dim}(F) = 4 - 2 = 2$, and are $S^1$-bundles over fibers of $f|_D$, which in turn are one-dimensional. The only such two-manifolds are unions of tori.
\ep
Recall that when $Z$ is separating it admits a global defining function $z \in C^\infty(Y)$. In this case, any boundary map $f$ to $(Y,Z)$ must map onto one component of $Y$ with respect to $Z$. More precisely, given such a $z$, consider $Y_+ = z^{-1}([0,\infty)) \subset Y$, a manifold with boundary given by $Z = z^{-1}(0)$. Then we have the following.
\begin{prop}\label{prop:globalbdrymap} Let $f\colon (X,D) \to (Y,Z)$ be a boundary map, and $Z$ separating. Then there exists a global defining function $z$ for $Z$ so that $f(X) \subset Y_+$, and $f$ defines a boundary map $f\colon (X,D) \to (Y', Z')$, where $Y' = Y_+ \cap f(X)$ and $Z' = Z \cap f(D)$.
\end{prop}
\bp Let $z$ be a defining function for $Z$. Then $f^*(z)$ is globally defined on $X$, with zero set $D = f^{-1}(Z)$. As $D$ has codimension at least two in $X$, its complement $X \setminus D$ is connected, so that $f^*(z)$ has fixed sign on $X \setminus D$. Consequently, by changing $z$ to $-z$ if necessary, the function $f^*(z)$ is non-negative. But then for all points $x \in X$, $z(f(x)) = f^*(z)(x) \geq 0$, so that $f(X) \subset Z_+$. Moreover, $f$ is a boundary map when restricting its codomain to its image.
\ep
\begin{rem}\label{rem:bdrymapname} The previous proposition explains the name `boundary map', as the defining condition specifies the behavior of $f$ normal to $Z$, the boundary of its restricted codomain.
\end{rem}
\subsection{Boundary (Lefschetz) fibrations}
We introduce further submersiveness assumptions.
\begin{defn} A \emph{boundary fibration} is a fibrating boundary map $f\colon (X,D) \to (Y,Z)$ such that $f|_{X \setminus D}\colon X \setminus D \to Y \setminus Z$ is a fibration.
\end{defn}
\begin{defn} A \emph{\placeholder} is a fibrating boundary map $f\colon (X^{2n},D) \to (\Sigma^{2},Z)$ such that $f|_{X \setminus D}\colon X \setminus D \to \Sigma \setminus Z$ is a Lefschetz fibration.
\end{defn}
The following is immediate from \autoref{cor:bdrymaplamorphism} together with \autoref{rem:fibrliealgdmorph}.
\begin{cor}\label{cor:bdryliealgdmorph} Let $f\colon (X,D) \to (Y,Z)$ be a codimension-two boundary (Lefschetz) fibration. Then $f$ induces a Lie algebroid (Lefschetz) fibration $(\varphi,f)\colon TX(-\log|D|) \to TY(-\log Z)$.
\end{cor}
The statement that $f$ induces a Lie algebroid morphism is to be interpreted as in \autoref{prop:elltologmorphism}, namely that there is a Lie algebroid morphism $(\varphi,f)$ such that $\varphi = df$ on sections. Moreover, the elliptic divisor structure on $D$ is the one induced from $f$ and $Z$.
By adapting the usual argument, we can ensure that boundary (Lefschetz) fibrations have connected fibers.
\begin{prop}\label{prop:blfconnfibers} Let $f\colon (X,D) \to (Y,Z)$ be a codimension-two boundary (Lefschetz) fibration, with $Z$ separating. Define $Y' := f(X)$ and $Z' := f(D)$. Then there exists a cover $g\colon \wt{Y}' \to Y'$ of $Y'$ and a boundary (Lefschetz) fibration $\wt{f}\colon (X,D) \to (\wt{Y}',Z')$ with connected generic fibers such that $f = g \circ \wt{f}$.
\end{prop}
\bp As $f$ is a global boundary map, by \autoref{prop:globalbdrymap} we obtain a boundary map $f\colon (X,D) \to (Y', Z')$ which immediately is also a boundary (Lefschetz) fibration. By definition, $f\colon X \setminus D \to Y' \setminus Z'$ is a (Lefschetz) fibration. Consequently, denoting its generic fiber by $F$, there is a sequence in homotopy $\pi_1(F) \to \pi_1(X \setminus D) \to \pi_1(Y' \setminus Z') \to \pi_0(F) \to 0$ \cite[Proposition 8.1.9]{GompfStipsicz99}. Applying Van Kampen's theorem for each connected component of $D$ shows that $\pi_1(X \setminus D)$ surjects onto $\pi_1(X)$ via the natural inclusion. Because $Z'$ is the boundary of $Y'$, we have $\pi_1(Y' \setminus Z') \cong \pi_1(Y')$. We obtain the following commutative diagram.
\begin{center}
	\begin{tikzpicture}
	\matrix (m) [matrix of math nodes, row sep=2.5em, column sep=2.5em,text height=1.5ex, text depth=0.25ex]
	{	 \pi_1(F) & \pi_1(X\setminus D) & \pi_1(Y' \setminus Z') & \pi_0(F) & 0\\ & \pi_1(X) & \pi_1(Y') & & \\};
	\path[-stealth]
	(m-1-1) edge node [left] {} (m-1-2)
	(m-1-2) edge node [left] {} (m-1-3)
	(m-1-3) edge node [left] {} (m-1-4)
	(m-1-4) edge node [left] {} (m-1-5)
	(m-2-2) edge node [above] {} (m-2-3)
	(m-1-3) edge node [left] {$\cong$} (m-2-3);
	\draw[>=stealth,->>]
	(m-1-2) edge node [left] {} (m-2-2);
	\end{tikzpicture}
\end{center}
The generic fiber $F$ of $f\colon X \setminus D \to Y' \setminus Z'$ is compact, hence $\pi_0(F)$ is finite. But then $f_*(\pi_1(X \setminus D)) \subset \pi_1(Y' \setminus Z')$ is a subgroup of finite index. Using the isomorphism $\pi_1(Y' \setminus Z') \cong \pi_1(Y')$, let $G$ denote the corresponding finite index subgroup inside $\pi_1(Y')$, and let $g\colon \wt{Y}' \to Y'$ be the associated cover. Then $f\colon X \to Y'$ lifts to $\wt{f}\colon X \to \wt{Y}'$ if and only if $f_*(\pi_1(X)) \subset G$, but this is an equality because $\pi_1(X \setminus D) \to \pi_1(X)$ surjects. For $\wt{f}$ we have that $\wt{f}\colon X \setminus D \to \wt{Y}' \setminus Z'$ induces a surjection $\wt{f}_*\colon \pi_1(X \setminus D) \to \pi_1(\wt{Y}' \setminus Z')$. But then $\pi_0(\wt{F})$ is trivial, so that the generic fibers $\wt{F}$ of $\wt{f}$ are connected.
\ep
By the above we can replace a boundary Lefschetz fibration on $(X,D)$ with $Z$ separating by one for which the generic fibers in $X \setminus D$ are connected. We next study the fibrating case, concluding that the fibers of $f|_D\colon D \to Z$ are connected if those of $f$ near $D$ are.
\begin{prop} Let $f\colon (X,D) \to (Y,Z)$ be a fibrating boundary map whose generic fibers near $D$ are connected, and $Z$ separating. Then the fibers of $f|_D\colon D \to Z$ are connected.
\end{prop}
\bp Using \autoref{prop:globalbdrymap}, replace $(Y,Z)$ by $(Y',Z')$ and let $V \subset Y'$ be an open neighbourhood of a point $y \in Z'$ so that $f|_{f^{-1}(V) \setminus D}$ is a fibration. Choose a curve $\gamma\colon [0,1] \to Y'$ such that $\gamma([0,1)) \in V \setminus Z'$, $f(1) = y$ and $\gamma$ is transverse to $Z'$. Then $M := f^{-1}(\gamma([0,1])) \subset f^{-1}(V)$ is a compact submanifold of $X$. Let $F$ be the generic fiber of $f$ near $D$ and denote $F_y = f^{-1}(y) = f|_D^{-1}(y)$. Let $F_{y,i}$ be the connected components of $F_y$, and let $U_i$ be disjoint opens around $F_{y,i}$ in $M$. Choose a sequence of points $(y_n)_{n \in \N} \subseteq \gamma([0,1])$ converging to $y$ such that $y_n \neq y$ for all $n \in \N$. For each $n$, consider the set $V_n := F_{y_n} \setminus \bigsqcup_i (U_i \cap F_{y_n})$. As each $U_i$ is open in $M$, the sets $U_i \cap F_{y_n}$ are disjoint opens in $F_{y_n}$.

Assume that $F_y$ is not connected. Then for $n$ large enough, there will be at least two indices $i$ for which $U_i \cap F_{y_n} \neq \emptyset$. As $F_{y_n}$ is connected, we conclude that $V_n \neq \emptyset$, as $F_n$ cannot be covered by disjoint opens. For each such $n$, let $x_n \in V_n$ be some element. By compactness of $M$, the sequence $\{x_n\}_{n \in \N}$ has a convergent subsequence, so that after relabeling we have $x_n \to x$ for some $x \in M$. By definition $x \in M \setminus (\bigcup_i U_i)$, so that $x \not\in U_i$ for all $i$, hence $x \not\in F_y$. However, by continuity we have $f(x_n) = y_n \to y = f(x)$, which is a contradiction. We conclude that $F_y$ is connected, so that the fibers of $f|_D\colon D \to Z'$ are connected.
\ep
\section{Constructing \placeholder{}s}
\label{sec:constrbdrlfs}
In this section we discuss how to obtain four-dimensional \placeholder{}s. In particular, we construct them out of genus one Lefschetz fibrations $f\colon X^4 \to \Sigma^2$ by surgery, replacing neighbourhoods of fibers of points $x \in \Sigma$ by a certain standard boundary map.

Given $n \in \Z$, let $p\colon L_n \to T^2$ be the complex line bundle over $T^2$ with first Chern class equal to $n \in H^2(T^2;\Z)$. Choose a Hermitian metric on $L_n$, which provides a fiberwise radial distance $r\colon {\rm Tot}(L_n) \to \R_{\geq 0}$. Further, let ${\rm pr}\colon T^2 \to S^1$ be the projection onto the first factor.
\begin{defn} The total space ${\rm Tot}(L_n)$ together with the map $f\colon{\rm Tot}(L_n) \to S^1 \times \R_{\geq 0}$ given by $f(x) := ({\rm pr}(p(x)), r^2(x))$ for $x \in {\rm Tot}(L_n)$, is called the \emph{standard $n$-model}.
\end{defn}
The first thing to note is that $f$ as defined above is an example of a boundary fibration.
\begin{prop}\label{prop:standardbdryfibration} The map $f\colon ({\rm Tot}(L_n), T^2) \to (S^1 \times \R_{\geq 0}, S^1 \times \{0\})$ is a boundary fibration, where $T^2$ is identified with the zero section in ${\rm Tot}(L_n)$.
\end{prop}
\bp Certainly $f$ is a strong map of pairs, and the codimension of $D = T^2$ inside ${\rm Tot}(L_n)$ is (at least) two. The normal Hessian of $H^\nu(f)$ of $f$ along $D$ is given by the constant map equal to $2$, which is clearly definite. We conclude that $f$ is a boundary map. Both the bundle projection $p$ as well as the projection ${\rm pr}$ are submersive, so that $f|_D$ is submersive, making $f$ a fibrating boundary map. It is immediate that $f$ is submersive when $r \neq 0$, i.e.\ on ${\rm Tot}(L_n) \setminus T^2$, so that $f$ is a boundary fibration.
\ep
Let $f\colon X \to Y$ be a smooth map that is a fibration in a neighbourhood of an embedded circle $\gamma \subset Y$. Recall that the \emph{monodromy} of $f$ around $\gamma$ is the element in the mapping class group of the fiber $F$ of $f^{-1}(\gamma) \to \gamma$ as a mapping torus.

Returning to the standard $n$-model, let $\varepsilon > 0$ be small and consider the circle $\gamma = S^1 \times \{\varepsilon\} \subset S^1 \times \R_{\geq 0}$. We compute the monodromy of $f\colon {\rm Tot}(L_n) \to S^1 \times \R_{\geq 0}$ around $\gamma$.
\begin{prop} \label{prop:monodromy}The monodromy of $f$ around $\gamma$ is the $-n$th power of a Dehn twist.
\end{prop}
\bp Let $M := f^{-1}(\gamma)$ be the $T^2$-bundle $f\colon M \to \gamma$. To compute the monodromy of $f$ around $\gamma$, note that using $p\colon M \to T^2$ we can view $M$ as the principal $\varepsilon$-sphere bundle of $L_n$. The Hermitian metric gives rise to a connection $i \theta \in \Omega^1(M;\R)$. As $L_n \to T^2$ has Euler class equal to $n$, its curvature equals $\frac{d\theta}{2 \pi i} = n da \wedge db$, where $da$ and $db$ are generators of $H^1(T^2)$. Recall now that circle bundles are classified by their Euler class due to the Gysin sequence. Consider $M' := \R^3 / \Gamma$, with $\Gamma = \langle \alpha_1,\alpha_2,\alpha_3\rangle$ the integral lattice generated by the following group elements acting on $\R^3$: 
\be
	\alpha_1(x,y,z) = (x,y+1,z), \qquad \alpha_2(x,y,z) = (x,y,z+1), \qquad \alpha_3(x,y,z) = (x+1,y, z-n y).
\ee
The projection ${\rm pr}_{12}\colon \R^3 \to \R^2$ given by $(x,y,z) \mapsto (x,y)$ descends to a map $g\colon \R^3 / \Gamma \to \R^2 / \Z^2$. This is an $S^1$-bundle over $T^2$, with invariant one-forms $e_1 = dx$ and $e_2 = dy$, and connection one-form $e_3 = dz + n x dy$. As $d e_3 = n e_1 \wedge e_2$, we conclude that $\R^3 / \Gamma$ has Euler class equal to $n$, so that $\R^3 / \Gamma \cong M$, as circle bundles over $T^2$. Note that the map ${\rm pr}_1\colon \R^3 \to \R$ given by $(x,y,z) \mapsto x$ also descends and exhibits $M'$ as a $T^2$-bundle over $\R / \Z \cong S^1$. To obtain a $T^2$-bundle over $S^1$ out of $p:\, M \to T^2$, the choice of any projection $T^2 \to S^1$ gives isomorphic bundles. Hence we can assume that the projection ${\rm pr}_1$ makes the following diagram commute.
\begin{center}
	\begin{tikzpicture}
	\matrix (m) [matrix of math nodes, row sep=2.5em, column sep=2.5em,text height=1.5ex, text depth=0.25ex]
	{	M & T^2 & S^1 \\ M' & T^2 & S^1 \\};
	\path[-stealth]
	(m-1-1) edge node [above] {$p$} (m-1-2)
	(m-1-2) edge node [above] {${\rm pr}$} (m-1-3)
	(m-2-1) edge node [above] {${\rm pr}_{12}$} (m-2-2)
	(m-2-2) edge node [above] {${\rm pr}_{1}$} (m-2-3)
	(m-1-1) edge node [left] {$\cong$} (m-2-1)
	(m-1-2) edge node [left] {$\cong$} (m-2-2)
	(m-1-3) edge node [left] {$\cong$} (m-2-3);
	\end{tikzpicture}
\end{center}
We can thus compute the monodromy of $f$ around $\gamma$ by considering ${\rm pr}_1:\, M' \to S^1$. Due to our concrete description it is immediate that
\be
M' = T^2 \times S^1 / (x,0) \sim (f(x),1), \qquad f \in {\rm Aut}(T^2) \text{ with } f_* = \begin{pmatrix} 1 & -n \\ 0 & 1 \end{pmatrix}\colon H_2(T^2) \to H_2(T^2).
\ee
We conclude that $M$ has monodromy equal to the $-n$th power of a Dehn twist as desired.
\ep
We next describe the surgery process whereby we replace suitable neighbourhoods of points by the above standard models.
\begin{defn} A \emph{punctured surface} is an open surface $\Sigma$ obtained from a closed surface by removing a finite number of disks.
\end{defn}
A punctured surface naturally has a compact closure $\overline{\Sigma}$ by adding the circle boundaries of the removed disks. This closure comes with a natural isomorphism $j\colon \Sigma \stackrel{\cong}{\to} \overline{\Sigma} \setminus \partial \overline{\Sigma}$.
\begin{prop}\label{prop:completion} Let $f\colon X^4 \to \Sigma^2$ be a proper map over a punctured surface such that $f$ is a genus one fibration in a neighbourhood of $\partial \Sigma$. Assume that the monodromy of $f$ around every boundary component is a power of a Dehn twist. Then there exists a compact elliptic pair $(\wt{X},|D|)$ with $i\colon X \stackrel{\cong}{\to} \wt{X} \setminus D$ and a fibrating boundary map $\wt{f}\colon (\wt{X},|D|) \to (\overline{\Sigma}, \partial \overline{\Sigma})$ such that the following diagram commutes:
	\begin{center}
		\begin{tikzpicture}
		\matrix (m) [matrix of math nodes, row sep=2.5em, column sep=2.5em,text height=1.5ex, text depth=0.25ex]
		{	X & (\wt{X},|D|) \\ \Sigma & (\overline{\Sigma}, \partial \overline{\Sigma}) \\};
		\path[-stealth]
		(m-1-1) edge node [left] {$f$} (m-2-1)
		(m-1-2) edge node [right] {$\wt{f}$} (m-2-2);
		
		\draw [right hook-latex]
		(m-1-1) edge node [above] {$i$} (m-1-2)
		(m-2-1) edge node [above] {$j$} (m-2-2);
		\end{tikzpicture}
	\end{center}
In the above situation, we say $f\colon X \to \Sigma$ \emph{can be completed} to $\wt{f}\colon (\wt{X},|D|) \to (\overline{\Sigma}, \partial \overline{\Sigma})$. These completions are unique when the monodromy is not trivial.
\end{prop}
\bp At each boundary end $E_i$ of $\Sigma$, glue in the standard $n$-model ${\rm Tot}(L_{n_i})$, where $n_i \in \Z$ is such that ${\rm Mon}(\gamma_i) = \delta^{n_i}$, with $\gamma_i = \partial E_i$. Here $\delta$ denotes a Dehn twist. As the monodromies agree, the fibrations are isomorphic, hence can be glued together to a new fibration.

It only remains to argue that if the monodromy is nontrivial, the completion is unique. Since the completion of each end of $X$ is obtained by gluing  $L_n$ to the end of $X$, by identifying $M$, the  $\varepsilon$-sphere bundle of $L_n$, with the pre-image of a loop around a boundary component of $\Sigma$, we see that any  other completion can be obtained by precomposing the gluing map by a diffeomorphism of $M$. However, due to a result of Waldhausen (\cite[Theorem 5.5]{Waldhausen67})  any such diffeomorphism  is isotopic to a fiber-preserving diffeomorphism (covering a diffeomorphism of the base) and therefore extends to the $\varepsilon$-disk bundle of $L_n$. Hence the completion of the end is unique up to diffeomorphism.
\ep
\begin{cor}\label{cor:lfcompletion} Let $f\colon X^4 \to \Sigma^2$ be a genus one \lf over a punctured surface. Assume that the monodromy of $f$ around every boundary component is a power of a Dehn twist. Then $f$ can be completed to a codimension-two \placeholder $\wt{f}\colon (\wt{X},|D|) \to (\overline{\Sigma}, \partial \overline{\Sigma})$.
\end{cor}
\bp The Lefschetz singularities lie in the interior of $\Sigma$ so that $f$ is a genus one fibration near $\partial \Sigma$. After completing using \autoref{prop:completion}, $\wt{f}$ is a fibrating boundary map and a genus one Lefschetz fibration in the interior, so that it is a \placeholder.
\ep
Consequently, given a genus one \lf on $X$, removing several disks and their inverse image and then completing as in \autoref{cor:lfcompletion} gives a \placeholder on $\wt{X}$. Homological essentialness of generic fibers is preserved by the completion process.
\begin{prop}\label{prop:retainhomess} Let $f\colon X^4 \to \Sigma^2$ be a genus one \lf with boundary Lefschetz completion $\wt{f}\colon (\wt{X},|D|) \to (\overline{\Sigma}, \partial \overline{\Sigma})$. Then $f$ is homologically essential if and only if $\wt{f}$ is.
\end{prop}
A codimension-two boundary \lf on a four-manifold $X$ has singular fibers equal to the Euler characteristic of $X$, as is true for genus one Lefschetz fibrations.
\begin{prop}\label{prop:blfsingfibers} Let $f:\, (X^4, D) \to (\Sigma^2, Z)$ be a codimension-two boundary Lefschetz fibration. Then we have $\chi(X) = \mu$, with $\mu$ the number of singular fibers of $f|_{X\setminus D}$.
\end{prop}
\bp Given two open sets $U, V \subseteq X$ we have $\chi(U \cup V) = \chi(U) + \chi(V) - \chi(U \cap V)$. Let $U := X \setminus D$ and take $V \subseteq ND$ to be an open subset in a tubular neighbourhood of $D$. Then $V$ is homotopy equivalent to $D$, which is a union of tori, so that $\chi(V) = \chi(D) = 0$. Moreover, $U \cap V$ is deformation equivalent to $S^1 ND$, which is a principal circle bundle over $D$, hence $\chi(U \cap V) = 0$ as well. Recall that if a manifold $M$ admits a genus $g$ Lefschetz fibration, its Euler characteristic satisfies $\chi(M) = 2 (2-2g) + \mu$, where $\mu$ is the number of singular fibers. Applying this to $U = X \setminus D$ with $g = 1$ we obtain that $\chi(X) = \chi(X \setminus D) = \mu$.
\ep
\section{Constructing \sgcs{}s}
\label{sec:constrsgcs}
In this section we combine the results from previous sections to construct \sgcs{}s using codimension-two boundary fibrations with two-dimensional fibers, and from codimension-two boundary \lf{}s in dimension four. For notational convenience, we call a boundary map $f\colon (X^{2n},D) \to (Y^{2n-2},Z)$ \emph{homologically essential} if the (generic) two-dimensional fibers $F$ (near $D$) specify nonzero homology classes ${[F] \neq 0 \in H_2(X\setminus D; \R)}$. We prove the following results stated in the introduction.
\begin{thm}\label{thm:lfibrationsgcs} Let $f\colon (X^4,D^2) \to (\Sigma^2,Z^1)$ be a homologically essential \placeholder between compact oriented manifolds such that $[Z] = 0 \in H_1(\Sigma;\Z_2)$. Then $(X,D)$ admits a \sgcs.
\end{thm}
\begin{thm}\label{thm:fibrationsgcs} Let $f\colon (X^{2n},D^{2n-2}) \to (Y^{2n-2},Z^{2n-3})$ be a homologically essential boundary fibration between compact oriented manifolds such that $(Y,Z)$ admits a log-symplectic structure. Then $(X,D)$ admits a \sgcs.
\end{thm}
In the remainder of this section, any mention of an elliptic divisor structure on $(X,D)$, or the Lie algebroid $TX(-\log |D|)$, will refer to the structure induced by a codimension-two boundary map via \autoref{prop:bdrymapmorphdivisors}. Similarly, we will without further mention have codimension-two boundary maps induce Lie algebroid morphisms between the appropriate elliptic and log-tangent bundles using \autoref{cor:bdryliealgdmorph}. Before proving the above two results, we first establish that we can create a suitable closed Lie algebroid two-form on $X$. 
\begin{prop}\label{prop:homessfibrbdrymap} Let $f\colon (X^{2n},D^{2n-2}) \to (Y^{2n-2},Z^{2n-3})$ be a homologically essential fibrating boundary map with $f(D)$ coorientable. Then there exists a closed two-form $\eta \in \Omega^2(\log |D|)$ with ${\rm Res}_q(\eta) = 0$ such that $\eta|_{\ker \varphi}$ is nondegenerate near $D$.
\end{prop}
\bp Note that $f$ is submersive with fibers $F$ around $D$ by \autoref{cor:fibratingsubm}. As $[F] \neq 0$, there exists $c \in H^2(X \setminus D)$ such that $\langle c, [F] \rangle = 1$. We construct local Lie algebroid two-forms around $D$ as in \autoref{prop:gluingfiberform}. By \autoref{thm:elllacohomology}, $H^2(X\setminus D)$ includes into $H_0^2(\log |D|)$. Let $\xi \in \Omega^2(\log |D|)$ be a closed two-form satisfying $[\xi] = c$ and ${\rm Res}_q(\xi) = 0$. By \autoref{prop:fibratingsemiglobal} there exist open neighbourhoods $U$ and $V$ around $D$ and $f(D)$ and a defining function for $f(D)$ on $V$ such that $f$ and $f|_D$ commute with the tubular neighbourhood projections ${\rm pr}_D$ and ${\rm pr}_Z$, and such that $f^*(z) = Q_{f^*(z)}$. Let $y \in f(D)$ and let $V' \subseteq V$ be a smaller open disk containing $y$ on which $NZ$ is trivial, and similarly so that $U' := f^{-1}(U)$. Then as in \autoref{prop:fibratingpointwise} we have $f(r e^{i\theta}, x, y) = (x, r^2)$. Moreover, in these coordinates we have $\ker \varphi = \langle \partial_\theta, \partial_y \rangle$. Let $\{V_i\}$ be a finite covering extracted from such open sets and let $U_i = f^{-1}(V_i)$, which together cover a neighbourhood of $D$. Set $U_0 = X \setminus D$. As $f\colon U_i \setminus D \to V_i \setminus Z$ is a $T^2$-bundle it is necessarily trivial, so that $H^2(U_i \setminus D) = H^2(T^2)$, with $H^1(T^2)$ generated by $\theta_i$ and $\gamma_i$ say.

Define $\eta_0 = \xi|_{X\setminus D}$ and $\eta_i \in \Omega^2(U_i;\log |D|)$ for $i \geq 1$ via $\eta_i = \lambda_i \theta_i \wedge \gamma_i$, where $\lambda_i \in \R$ is chosen such that $\int_F \eta_i = \int_F c$. Then $[\eta_i] = c|_{U_i} = [\xi|_{U_i}] \in H^2(U_i; \log |D|)$, so that there exist one-forms $\alpha_i$ such that $\eta_i = \xi + d \alpha_i$. As in \autoref{prop:gluingfiberform}, define a closed Lie algebroid two-form $\eta:= \xi + d (\sum_i (\psi_i\circ f) \alpha_{y_i})$ using a partition of unity $\{\psi_i\}$ subordinate to $\{V_i\}$. As ${\rm Res}_r(\eta_i) = 0$, we conclude immediately that ${\rm Res}_q(\eta_i) = 0$. Moreover, as each form $\alpha_{y_i}$ is smooth, this implies that ${\rm Res}_q(\eta) = 0$ as well. Finally, near $D$ the form $\eta$ is nondegenerate on $\ker \varphi$ as there it is given by the convex combination of forms $\eta = \sum_i (\psi \circ f) \eta_i$, where each $\eta_i$ is nondegenerate on $\ker \varphi = \langle \partial_\theta, \partial_y \rangle$ by construction.
\ep
We can now prove \autoref{thm:lfibrationsgcs} and \autoref{thm:fibrationsgcs} simultaneously, using \autoref{thm:thurstonfibration} or \autoref{thm:thurstonlfibration} respectively.
\bp[ of \autoref{thm:lfibrationsgcs} and \autoref{thm:fibrationsgcs}] Let $\omega_Y$ be a log-symplectic structure on $(Y,Z)$. As $Y$ is oriented, by \autoref{prop:blogbsymp} and \autoref{prop:orlogseparating} we know that $Z$ is separating, and in particular coorientable. Using \autoref{prop:globalbdrymap}, replace $(Y,Z)$ by $(Y',Z')$, and then using \autoref{prop:blfconnfibers}, lift to a cover $g\colon \wt{Y}' \to Y'$ so that $\wt{f}\colon (X,D) \to (\wt{Y}',Z')$ has connected fibers. It is immediate that $g^*\omega_Y$ defines a log-symplectic structure on $(\wt{Y}',Z)$, and moreover that $\wt{f}$ is a boundary Lefschetz fibration. The generic fibers of  $\wt{f}\colon X \to \wt{Y}'$ are either all homologically essential, or none are. Let $y \in Y' \setminus Z'$ be a regular point and let $y_1, \dots, y_n \in g^{-1}(y) \subset \wt{Y}'$ be its inverse images under the covering, with $F_{y_1}, \dots F_{y_n}$ their fibers under $\wt{f}$. Then for some choice of signs we have $[F_{y_1}] \pm \dots \pm [F_{y_n}] = [F]$. However, $[F] \neq 0$ by assumption, so there exists $i \in \{1,\dots,n\}$ such that $[F_{y_i}] \neq 0$. But then $[F_{y_i}] \neq 0$ for all $i \in \{1,\dots,n\}$, so $\wt{f}$ is homologically essential.

Apply \autoref{prop:homessfibrbdrymap} to the boundary (Lefschetz) fibration $\wt{f}$. This yields a global closed Lie algebroid two-form $\eta \in \Omega^2_{\rm cl}(X, \log |D|)$ such that $\eta|_{\ker \varphi}$ is nondegenerate near $D$. Moreover, its cohomology class $c = [\eta]$ pairs nontrivally with (generic) fibers in $X \setminus D$. Now apply either \autoref{thm:thurstonfibration} or \autoref{thm:thurstonlfibration}, possibly changing $\eta$ to a form $\eta'$ which agrees with $\eta$ in a neighbourhood of $D$. We obtain an elliptic symplectic structure $\omega_t = \varphi^* \omega_Y + t \eta'$ on $(X,|D|)$. As ${\rm Res}_q(\eta) = 0$ and $\eta'$ agrees with $\eta$ near $D$, we have ${\rm Res}_q(\eta') = 0$ as well. Using \autoref{prop:residuemaps} we see that ${\rm Res}_q(\omega_t) = {\rm Res}_q(\varphi^*\omega_Y) + t\, {\rm Res}_q(\eta') = 0$. But then the elliptic symplectic structure $\omega_t$ has zero elliptic residue. By \autoref{thm:sgcscorrespondence}, we conclude that $\omega_t$ for $t > 0$ small determines a \sgcs on $(X,|D|)$.
\ep
The \sgcs{}s on $(X,D)$ constructed in the proofs of \autoref{thm:lfibrationsgcs} and \autoref{thm:fibrationsgcs} arise as elliptic log-symplectic forms $\omega_X$ with zero elliptic residue through the correspondence of \autoref{thm:sgcscorrespondence}. As mentioned below \autoref{thm:sgcscorrespondence}, the three-form $H$ required for integrability of the \gcs is given by $H = {\rm Res}_r(\omega_X) \cup {\rm PD}_X[D]$. The two-form $\eta'$ introduced during the proof satisfies ${\rm Res}_r(\eta') = 0$. Using \autoref{prop:residuemaps}, together with the fact that $\omega_X = \varphi^* \omega_Y + t \eta'$ for some $t > 0$, we see that ${\rm Res}_r(\omega_X) = {\rm Res}_r(f^* \omega_Y) + t\, {\rm Res}_r(\eta') = f^* {\rm Res}_Z(\omega_Y)$.
\begin{rem} \label{rem:7.4}Given a defining function $z$ for $Z$ we can locally write $\omega_Y = d \log z \wedge \alpha + \beta$, with $(\alpha,\beta) = ({\rm Res}_Z(\omega_Y), j_Z^* \omega_Y)$ the induced cosymplectic structure on $Z$. Then $f^*\omega_Y = d\log r \wedge f^*(\alpha) + f^*(\beta)$, and ${\rm Res}_r(f^*\omega_Y) = f^*({\rm Res}_Z(\omega_Y)) = f^*(\alpha)$.
\end{rem}
\subsection{Fibering over $T^2$ and $S^1$}
Recall from Section \ref{sec:gcsstablegcs} that in the compact case, for a stable generalized pair $(X,D)$ and a log-Poisson pair $(Y,Z)$, the manifolds $D$ and $Z$ fiber over $T^2$ and $S^1$ respectively. After perturbing each structure, this can be achieved using one-forms on $D$ and $Z$ induced by the geometric structures. For a \sgcs $\mc{J}$ on X, these one-forms are $({\rm Res}_r,{\rm Res}_\theta)(\omega_X) \in \Omega^{1}(D) \times \Omega^1(D)$, where $\omega_X$ is the elliptic symplectic structure obtained from $\mc{J}$ using \autoref{thm:sgcscorrespondence}, see \cite{CavalcantiGualtieri15}. In the log-Poisson case, the one-form used instead is ${\rm Res}_Z(\omega_Y)$ (see \cite{GualtieriLi14, GuilleminMirandaPires14, MarcutOsornoTorres14, OsornoTorres15}), with $\omega_Y$ the log-symplectic structure induced by the log-Poisson structure using \autoref{prop:blogbsymp}.
\begin{rem} When $Z$ fibers over $S^1$ using $\alpha := {\rm Res}_Z(\omega_Y)$, the form $\alpha$ is Poincar\'e dual on $Z$ to the fiber $F_p$ of the induced fibration $p_Z\colon Z \to S^1$. In this case we obtain that $H = f^*({\rm PD}_Z[F_p]) \cup {\rm PD}_X[D]$.
\end{rem}
The \sgcs{}s we construct out of boundary (Lefschetz) fibrations $f\colon (X,D) \to (Y,Z)$  are such that $D$ fibers over $Z$. It is natural to ask whether this can be made compatible with the above fibrations $p_D\colon D \to T^2$ and $p_Z\colon Z \to S^1$. Note that $({\rm Res}_r,{\rm Res}_\theta)(\omega_X) = (f^*{\rm Res}_Z(\omega_Y), t\, {\rm Res}_\theta(\eta))$ as $\omega_X = \varphi^* \omega_Y + t \eta$, and ${\rm Res}_\theta(\varphi^* \omega_Y) = 0$ by the discussion below \autoref{prop:bdrymaplocalcoord}. This immediately shows we can make the fibrations compatible, as we can choose $\eta$ such that $t\, {\rm Res}_\theta(\eta)$  and $f^*{\rm Res}_Z(\omega_Y)$ determine a fibration, instead of a just foliation. As in Section \ref{sec:constrbdrlfs}, let $p\colon T^2 \to S^1$ be projection onto the first factor.
\begin{cor}\label{cor:compatiblefibrations} Under the assumptions of \autoref{thm:lfibrationsgcs} or \autoref{thm:fibrationsgcs}, given a \blog $\omega_Y$ on $(Y,Z)$, the elliptic two-form $\omega_X$ can be chosen such that $p \circ p_D = p_Z \circ f|_D$.
\end{cor}
In other words, we have a full commutative diagram around $D$ and $Z$ as follows:
\begin{center}
	\begin{tikzpicture}
	\matrix (m) [matrix of math nodes, row sep=2.5em, column sep=2.5em,text height=1.5ex, text depth=0.25ex]
	{	 ND & D & T^2 \\ NZ & Z & S^1\\};
	\path[-stealth]
	(m-1-1) edge node [above] {${\rm pr}_D$} (m-1-2)
	(m-1-1) edge node [left] {$f$} (m-2-1)
	(m-1-2) edge node [above] {$p_D$} (m-1-3)
	(m-1-2) edge node [left] {$f|_D$} (m-2-2)
	(m-2-1) edge node [above] {${\rm pr}_Z$} (m-2-2)
	(m-1-3) edge node [left] {$p$} (m-2-3)
	(m-2-2) edge node [above] {$p_Z$} (m-2-3);
	\end{tikzpicture}
\end{center}
\section{Examples and applications}
\label{sec:applications}
In this section we discuss some applications of the results obtained in this paper. We start with the following immediate consequence of the results in Sections \ref{sec:constrbdrlfs} and \ref{sec:constrsgcs}.
\begin{thm}\label{thm:g1lfsgcs} Let $f\colon X \to D^2$ be a genus one Lefschetz fibration over the disk whose monodromy around the boundary is a power of a Dehn twist. Then all possible completions $\wt{f}\colon (\wt{X}, D) \to (D^2, \partial D^2)$ admit a \sgcs.
\end{thm}
\bp Any Lefschetz fibration over the disk $D^2$ is homologically essential. By the monodromy assumption around the boundary, $f$ admits completions $\wt{f}$ using \autoref{cor:lfcompletion}. By \autoref{prop:retainhomess}, the map $\wt{f}$ is a homologically essential boundary Lefschetz fibration, so that $(\wt{X}, D)$ admits a \sgcs by \autoref{thm:lfibrationsgcs}.
\ep
We now turn to exhibiting \sgcs{}s on concrete four-manifolds using boundary (Lefschetz) fibrations.
\begin{exa} Let $I = [-1,1]$ with coordinate $t$ and view $S^2 \subseteq \R^3$ with north and south pole $p_N$ and $p_S$. Consider the standard height function $p\colon (S^2, p_N \cup p_S) \to (I,\partial I)$ given by $(x,y,z) \mapsto z$. The map $h\colon (I,\partial I) \to (\R, 0)$ given by $h(t) := 1 - t^2$ is a defining function for $\partial I$. Then $p^*(h)(x,y,z) = 1 - (1 - x^2 - y^2) = x^2 + y^2$, so that $p$ is a boundary map by the description of \autoref{prop:bdrymaplocalcoord}. In fact, $p$ is a boundary fibration. Let $\varphi\colon S^3 \to S^2$ be the Hopf fibration, making $p \circ \varphi\colon (S^3, S^1_N \cup S^1_S) \to (I,\partial I)$ a boundary fibration, where $S^1_N = \varphi^{-1}(p_N)$ and similarly for $S^1_S$. Finally, define the boundary fibration $f\colon (S^3 \times S^1, D) \to (I \times S^1, Z)$, where $Z = \partial I \times S^1$ and $D = (S^1_N \cup S^1_S) \times S^1 = f^{-1}(Z)$, by $f(x,\theta) = (p \circ \varphi)(x)$. Using \autoref{lem:blogsurface}, $(I \times S^1, Z)$ admits a \blog, while the fibers of $f$ are clearly homologically essential. By \autoref{thm:fibrationsgcs} applied to $f$, we conclude that $(S^1 \times S^3, D)$ admits a \sgcs\ whose type change locus has two connected components. This structure is integrable with respect to the zero three-form due to \autoref{rem:7.4}.
\end{exa}
Using a slightly different map, we can ensure that the type change locus is connected.
\begin{exa} View $S^3 \subseteq \C^2$ with coordinates $(z_1,z_2)$ in the standard way and consider the map $\varphi\colon S^3 \to D^2$ given by $\varphi(z_1,z_2) = z_2$, viewing $D^2 \subseteq \C$. In the interior $D^2 \setminus \partial D^2$, the map $\varphi$ admits a section $s(z) = (\sqrt{1 - |z|^2},z)$. As $\varphi \circ s = {\rm id}$, we see that $\varphi$ is submersive on $S^3 \setminus D$, where $D$ is the circle $D = \varphi^{-1}(\partial D^2) = \{(0,z_2) \in \C^2: |z_2| = 1\}$. Moreover, $\varphi|_D$ is the identity from $D$ to $Z = \partial D^2$, which in particular implies that $\varphi|_D$ is submersive. Let $h(z) := 1 - |z|^2$ be a defining function for $\partial D^2$ on $D^2$. Then $\varphi^*(h)(z_1,z_2) = 1 - |z_2|^2$, which shows that $\varphi$ is a boundary map using \autoref{prop:bdrymaplocalcoord}. But then $\varphi$ is in fact a boundary fibration. Noting that $\varphi|_{S^3 \setminus D}$ admits a section, $\varphi:\, S^3 \setminus D \to D^2 \setminus D^2$ is a trivial $S^1$-bundle. Now define a map of pairs $f\colon (S^3\times S^1, D \times S^1 ) \to (D^2, \partial D^2)$ by $f(z_1,z_2,x) := \varphi(z_1,z_2)$. This is also a boundary fibration, which is homologically essential as $f|_{ (S^3 \setminus D) \times S^1}$ is the projection $S^1  \times (D^2 \setminus \partial D^2) \times S^1 \to (D^2, \partial D^2)$. Finally, $(D^2, \partial D^2)$ admits a \blog by \autoref{lem:blogsurface}. Using \autoref{thm:fibrationsgcs} we conclude that $(S^3 \times S^1,  D \times S^1)$ admits a \sgcs with connected type change locus which is  integrable with respect to a nonzero degree three cohomology class.
\end{exa}
Before we consider actual boundary Lefschetz fibrations, it is convenient to reinterpret \autoref{prop:completion} with a Kirby calculus point of view. There, the completion of each end of $X$ was done by gluing in a copy of the total space of $L_n$, the disk bundle over the torus with Euler class $n$. Observe that ${\rm Tot}(L_n)$ is a four-manifold built out of one 0-handle, two 1-handles and one 2-handle. The process of capping off an end of $X$ with ${\rm Tot}(L_n)$ amounts to inverting the handle structure of ${\rm Tot}(L_n)$ and adding it to $X$. That is, to cap each end of $X$ we add one 2-handle, two 3-handles and one 4-handle. Further, due to inversion, the 2-handle of ${\rm Tot}(L_n)$ has its original core and co-core interchanged and once this 2-handle is added there is a unique way to complete with the 3- and 4-handles, up to diffeomorphism. Therefore we only have to describe how to attach the 2-handle, whose core is the circle fiber of the projection $M \to T^2$, where $M$ is the $\varepsilon$-sphere bundle of $L_n$. From the discussion in \autoref{prop:monodromy}, this fiber is a circle left fixed by the monodromy map, which, if the monodromy is nontrivial, is determined by the monodromy map and is essentially unique. Further, since the circle is the fiber of the projection $M \to T^2$, the framing, in double strand notation, is  given by the nearby fibers.
\begin{exa} Consider the genus one Lefschetz fibration over $D^2$ with one singular fiber with vanishing cycle $b$, where $b \in H^1(T^2)$ is a generator. Then the monodromy around $\partial D^2$ is $b$, which is clearly the power of a Dehn twist. The resulting completion $\wt{X}$ thus admits a \sgcs by \autoref{thm:g1lfsgcs}, and we are left with determining its diffeomorphism type. Draw a Kirby diagram for the trivial $T^2$-bundle over $D^2$ in the standard way (see \cite{GompfStipsicz99}) and attach a $-1$-framed two-handle to represent the vanishing cycle $b$. The completion process adds a single $0$-framed two-handle along the same cycle (see Figure \ref{fig:S1xS3}). Slide the $-1$-framed two-handle over this to split off a $-1$-framed unknot, representing a copy of $\overline{\C P}^2$. The remaining diagram collapses to a single one-handle after two further handle slides, which shows that $\wt{X} = S^1 \times S^3 \# \overline{\C P}^2$.
\begin{figure}[h!!]
\begin{center}
\includegraphics[height=4.7cm]{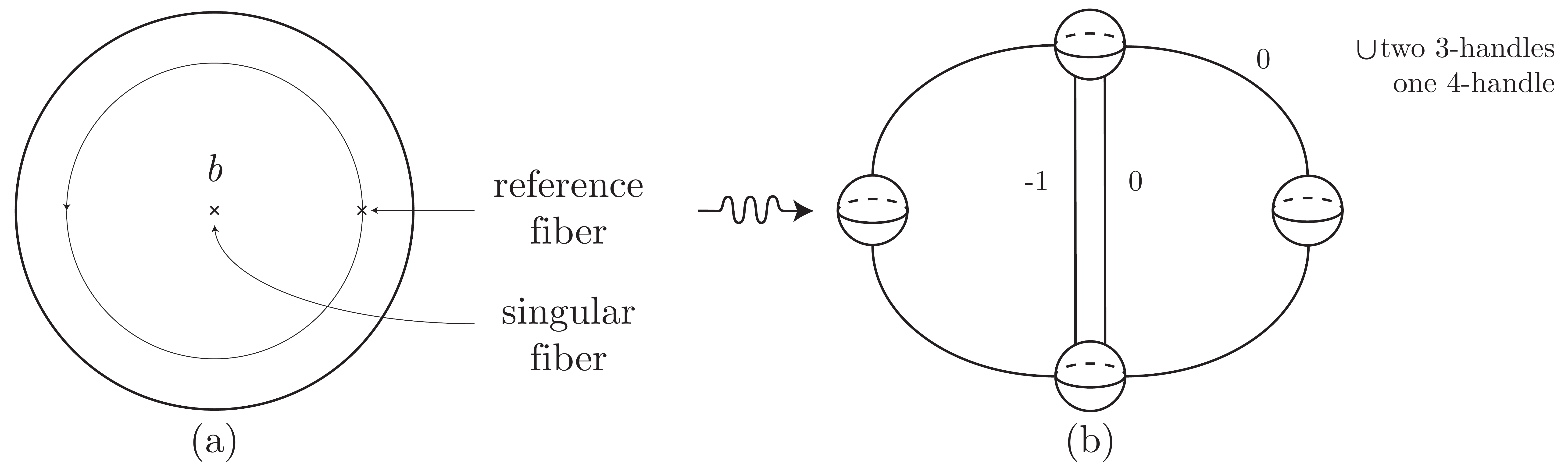} \caption{Figure (a) shows the base of a boundary Lefschetz fibration over the disk with a single nodal fiber, a choice of regular fiber for reference and loop around the boundary used to compute the monodromy and (b) shows the corresponding Kirby diagram.} \label{fig:S1xS3}
\end{center}
\end{figure}
\end{exa}
We can also recover the examples of \cite{CavalcantiGualtieri09}, showing that $m \C P^2 \# n \overline{\C P}^2$ admits a \sgcs if and only if $m$ is odd (i.e., if it admits an almost-complex structure).
\begin{exa} Let $m, n \in \N$ and assume that $m = 2k + 1$ is odd. Let $a, b \in H^1(T^2)$ be generators and consider $(D^2, \partial D^2)$ with clockwise assignment of singular fibers given by $a + (4k-1)b$, $a + 4 j b$ for $j$ from  $k-1$ up to $1-k$, and $a - (4k-1)b$, and $n$ copies of $b$. The monodromy around the boundary can be computed to be $(10k-1-n)b$. Therefore the associated genus one Lefschetz fibration admits a completion to $f\colon (\wt{X},D) \to (D^2, \partial D^2)$ (see Figure \ref{fig:An}). As is shown in \cite{CavalcantiGualtieri09}, we have $\wt{X} = m\C P^2\# n\overline{\C P}^2$, which hence has a \sgcs.
\begin{figure}[h!!]
\begin{center}
\includegraphics[height=5.5cm]{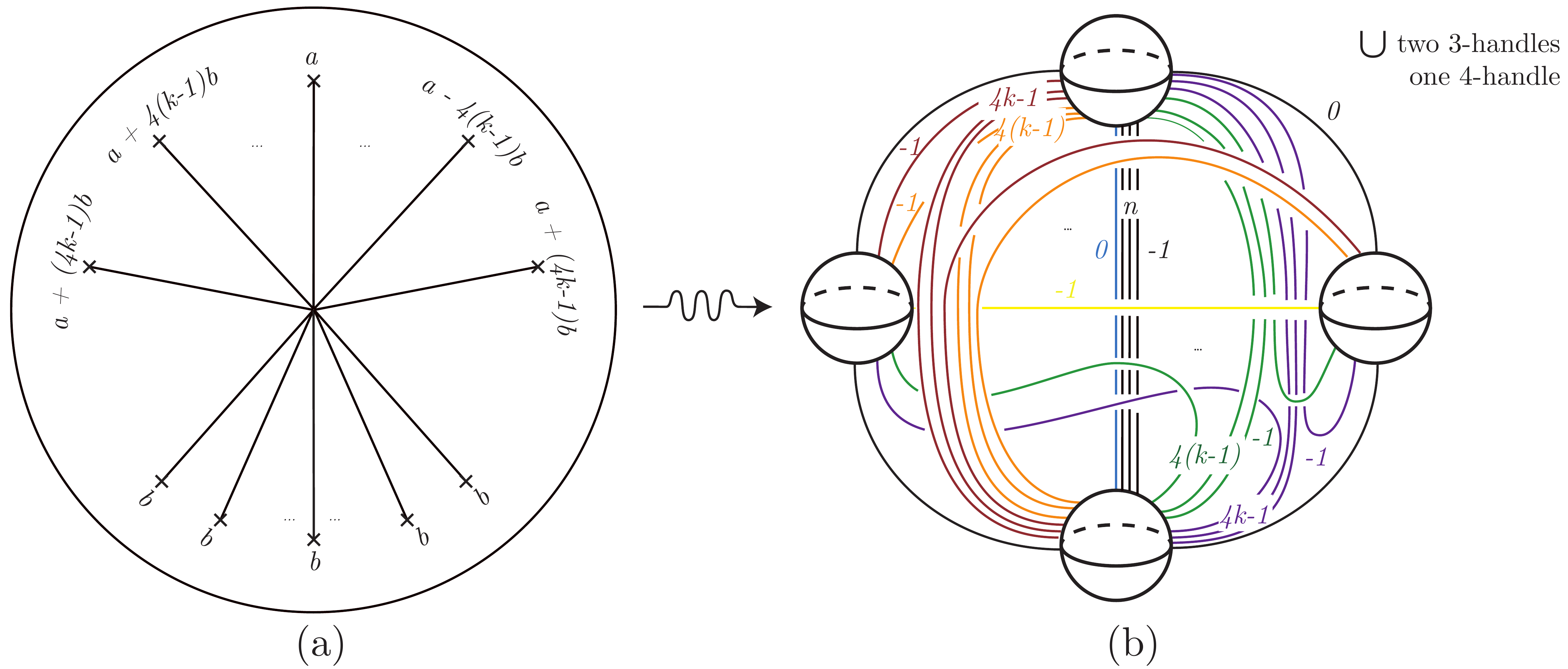} \caption{Figure (a) shows the base of a boundary Lefschetz fibration over the disc with $m+ n$ singular fibers and a choice of regular fiber for reference at the center of disk. Figure (b) shows the corresponding Kirby diagram.} \label{fig:An}
\end{center}
\end{figure}
\end{exa}
Note that not every \sgcs comes from a boundary \lf, similarly to the case for symplectic structures and Lefschetz fibrations, or \blog{}s and achiral Lefschetz fibrations (see \cite{CavalcantiKlaasse16}).
\begin{exa} The symplectic manifold $(\C P^2 , \omega_{\rm FS})$ carries a \sgcs with $D = \emptyset$. As $(\C P^2, \omega_{\rm FS})$ is not a symplectic fibration over any surface, this \sgcs can not be obtained through our construction.
\end{exa}
\begin{exa} Let $\Sigma_g$ be the genus $g$ surface and $g > 1$. Then $X = S^2 \times \Sigma_g$ has negative Euler characteristic so cannot admit a \placeholder by \autoref{prop:blfsingfibers}. Hence $X$ does not admit a \sgcs obtained through our methods. However, $X$ is symplectic hence carries a \sgcs with $D = \emptyset$.
\end{exa}
%

\bibliographystyle{hyperamsplain-nodash}
\bibliography{stablegcs}

\end{document}